\documentclass[letterpaper, 10 pt, journal, twoside]{ieeetran}

\IEEEoverridecommandlockouts                              % This command is only needed if
\usepackage{graphicx}
\graphicspath{ {figures/} }
\usepackage{subcaption}
\usepackage{lipsum}

\usepackage{color,soul}

\usepackage{amsmath}
\usepackage{amsfonts}
\usepackage{multirow}
\usepackage{multicol}
\usepackage[ruled]{algorithm2e}
\SetKwInput{KwInput}{Input}                % Set the Input
\SetKwInput{KwOutput}{Output}              % set the Output
\makeatletter
\let\NAT@parse\undefined
\makeatother
\usepackage{hyperref}
\newcommand{\Resp}[1]{{\color{black}  #1}}

\title{Decentralized Trajectory Optimization \\for Multi-Agent Ergodic Exploration
}

\author{Dimitris Gkouletsos, Andrea Iannelli, Mathias Hudoba de Badyn, John Lygeros% <-this % stops a space
\\\medskip\vspace{0.1cm}
\fbox{
Accepted for publication: Robotics and Automation Letters. (DOI: 10.1109/LRA.2021.3094242) \textcopyright 2021 IEEE.
}
\vspace{-0.8cm}
\thanks{Manuscript received: February, 24, 2021; Revised May, 26, 2021; Accepted June, 22, 2021.}%Use only for final RAL version
\thanks{This paper was recommended for publication by Editor Lucia Pallottino upon evaluation of the Associate Editor and Reviewers' comments.
This work has been partially supported
by the Swiss National Science Foundation under NCCR Automation.
}% <-this % stops a space
\thanks{
The authors are with the Department of Information Technology and Electrical Engineering, Automatic Control Lab, ETH, Z\"{u}rich 8092, Switzerland. (\emph{corresponding author: Andrea Iannelli})
emails: {\tt\small dgkouletsos@student.ethz.ch, \tt\small \{iannelli, mbadyn, lygeros\} @control.ee.ethz.ch}
}%
\thanks{Digital Object Identifier (DOI): see top of this page.}
}

\begin{document}

\maketitle
%\thispagestyle{empty}
%\pagestyle{empty}

%%%%%%%%%%%%%%%%%%%%%%%%%%%%%%%%%%%%%%%%%%%%%%%%%%%%%%%%%%%%%%%%%%%%%%%%%%%%%%%%
\begin{abstract}
Autonomous exploration is an application of growing importance in robotics. A promising strategy is ergodic trajectory planning, whereby an agent spends in each area a fraction of time which is proportional to its probability information density function. In this paper, a decentralized ergodic multi-agent trajectory planning algorithm featuring limited communication constraints is proposed. The agents' trajectories are designed by optimizing a weighted cost encompassing ergodicity, control energy and close-distance operation objectives. To solve the underlying optimal control problem, a second-order descent iterative method coupled with a projection operator in the form of an optimal feedback controller is used. Exhaustive numerical analyses show that the multi-agent solution allows a much more efficient exploration in terms of completion task time and control energy distribution by leveraging collaboration among agents.
%To set up the decentralized network topology, Erdos-Renyi random graphs are used. The agents exploration trajectories are derived through an optimization problem that investigates simultaneously  ergodicity, control energy and the closeness of the agents. To solve the optimal problem, we use an iterative descent method coupled with a projection operator in the form of an optimal feedback controller.  To enhance the benefits of multiple agents against single-agent systems, an
% exhaustive Monte Carlo analysis is performed. It is shown that multiple agents outperform single agents and collaborate to explore efficiently an area in terms of completion task time, control energy and traveled distance.

\end{abstract}

\section{INTRODUCTION}
\IEEEPARstart{I}{n} recent years, autonomous exploration has received significant attention in view of a variety of application fields such as agriculture surveillance
\cite{agriculture} and active map searching for disaster enviroments \cite{delmerico2017active}. %and underwater facility monitoring \cite{hollinger2013active}.
A promising framework in which to set the problem is ergodic trajectory planning \cite{Mathew}, whereby the generated trajectory samples a region in the search space proportional to the expectation of how informative the region will be.
% an ergodic metric defined as the deviation between a spatial information density function and a time-averaged
%systems trajectory distribution is introduced.
Intuitively, an ideal ergodic trajectory should cover high-valued information regions proportionally to the time spent in that region. A projection-based iterative optimization algorithm was proposed in \cite{Miller} to plan ergodic trajectories of a single-agent system. Experimental results \cite{miller2015ergodic} showcased the advantages of this approach, which outperformed alternative entropy minimization and information maximization strategies. In \cite{ayvali2017ergodic}, an extension was presented for constrained environments and obstacle avoidance problems, and the Kullback-Leibler divergence was adopted as alternative to the ergodic metric.

Recent advances in computational resources and maturity of distributed algorithms have enabled multi-agent network collaborations in a variety of critical applications such as object secure transportation \cite{lee_planning} and search-and-prosecute missions \cite{manathara}. Multi-agent configurations have been shown to outperform single-agent ones in terms of completion mission time \cite{manathara,voronoi}.
%In particular, decentralized solutions are of great appeal, since the presence of a central unit driving the decisions of the individuals is replaced by a network of agents planning in parallel the task execution \cite{bakule}.
Decentralized solutions are of great appeal, since replacing a central unit by a network of agents planning in parallel the task execution %\cite{bakule}.
%This
gives more robustness against failures and provides more communication flexibility \cite{nestmeyer2017decentralized}.

In this paper, we present a multi-agent decentralized trajectory planning problem in the framework of ergodic exploration. Related works are: \cite{abraham2018decentralized}, which reformulates the single-agent problem in \cite{Miller} using a Nash equillibrium interpretation for the multi-agent setting; \cite{salman2017multi}, focusing on area coverage with obstacles; \Resp{and \cite{prabhakar} which demonstrates a decentralized ergodic swarm control framework adaptable to external user commands and dynamic environmental information}.  \. %Further differences with the present work are highlighted in the following.
The contribution of this work is twofold. First, a decentralized multi-agent extension of the approach in \cite{Miller} is proposed (Section \ref{sIII}). The algorithm consists of four steps all performed at agent-level: a second-order steepest descent optimizer, combined with a line-search scheme, determines a candidate trajectory that is optimal according to a generalized global cost function (discussed in \ref{sIII}); %determines agents optimal trajectories to explore ergodically a domain. %A  is applied such that the agents satisfy their
a feasible trajectory that satisfies the dynamic constraints of the agent is obtained by projection \cite{Hauser}; estimates of the other agents trajectories are updated by averaging.
%Satisfaction of the individual dynamics constraints is achieved using the projection-based strategy from \cite{Hauser}.
In contrast to previous works, the cost determining the optimal trajectory encompasses: a (global) ergodic metric; a (global) control energy index; and a penalty on the inter-agent distance. % To the best of our knowledge, it is the first time that a distance cost argument is included that penalizes closeness of agents in the context of ergodic exploration.
\Resp{This new cost definition, and the applicability of the approach to general nonlinear systems without requiring control affine dynamics, are important differences with respect to \cite{abraham2018decentralized}, where the only objective was to optimize ergodicity and thus could make use of a different algorithmic approach.}
%mode insertion gradient algorithm was used to determine trajectories achieving a sufficient ergodic metric reduction.
%Moreover, our method applies to general nonlinear dynamical systems without any assumption on control affine dynamics
%\cite{abraham2018decentralized} or linear dynamics \cite{salman2017multi}.
Second, differently from previous works, we provide a systematic investigation of the benefits of multi-agent systems against single-agent systems in the ergodic exploration framework (Section \ref{sV}). Even though experimental results are not provided here, it is believed that the reported numerical analyses and the insights gained therein represent an invaluable starting point to plan a future experimental study, similar in concept to the one reported in \cite{miller2015ergodic} for the single-agent case.
%In spite of the significant research efforts in multi-agent ergodic exploration \cite{abraham2018decentralized}, \cite{salman2017multi},
%Hence, we employ an exhaustive Monte Carlo simulation to compare multi-agent and single-agent exploration  over
A range of performance metrics, including completion task time and control energy, are reported for different numbers of agents and network topologies randomly generated. \Resp{An empirical study is also performed to describe the convergence property of the algorithm}. %, defined using
%To develop a fully decentralized architecture, we define a network topology according
%
%Erdos-Renyi random graphs \cite{erdos}.
%Section \ref{sIV}, which details all the practical aspects of the algorithm for the considered scenarios, accompany the release of the \href{https://gitlab.ethz.ch/andreaia/multi_agent_explore}{repository}\footnote{{\scriptsize \texttt{https://gitlab.ethz.ch/andreaia/multi\_agent\_explore}}} that reproduces the results.
Section \ref{sIV} details all the practical aspects of the algorithm for the considered scenarios and accompanies the release of the repository \cite{repo_ETH} that reproduces the results.

%in the repository \url{https://gitlab.ethz.ch/andreaia/multi_agent_explore}.

%The use of random graphs indicates independence  of the algorithm  with
%respect to the graph topology instead of spatially connected agents as presented in other decentralized explorations
%\cite{nestmeyer2017decentralized}.

%The outline of this work is organized as follows: Section \ref{sII} introduces the ergodicity notion and Section \ref{sIII} analyzes the main steps of the
%algorithm. In Section \ref{sIV}, we present the implementation details employed in this work  to proceed with numerical simulations.
%exhibits the case studies of the proposed algorithm and discussion about the results, while Section \ref{sVI} is the conclusion.

%%%%%%%%%%%%%%%%%%%%%%%%%%%%%%%%%%%%%%%%%%%%%%%%%%%%%%%%%%%%%%%%%%%%%%%%%%%%%%%%

%ERGODICITY

\section{ERGODICITY}\label{sII}
\begin{figure}[h!] % [b!]
\centering
%\framebox{\includegraphics[scale=0.25]{Ergodicity_concept}}
\includegraphics[width=0.6\columnwidth]{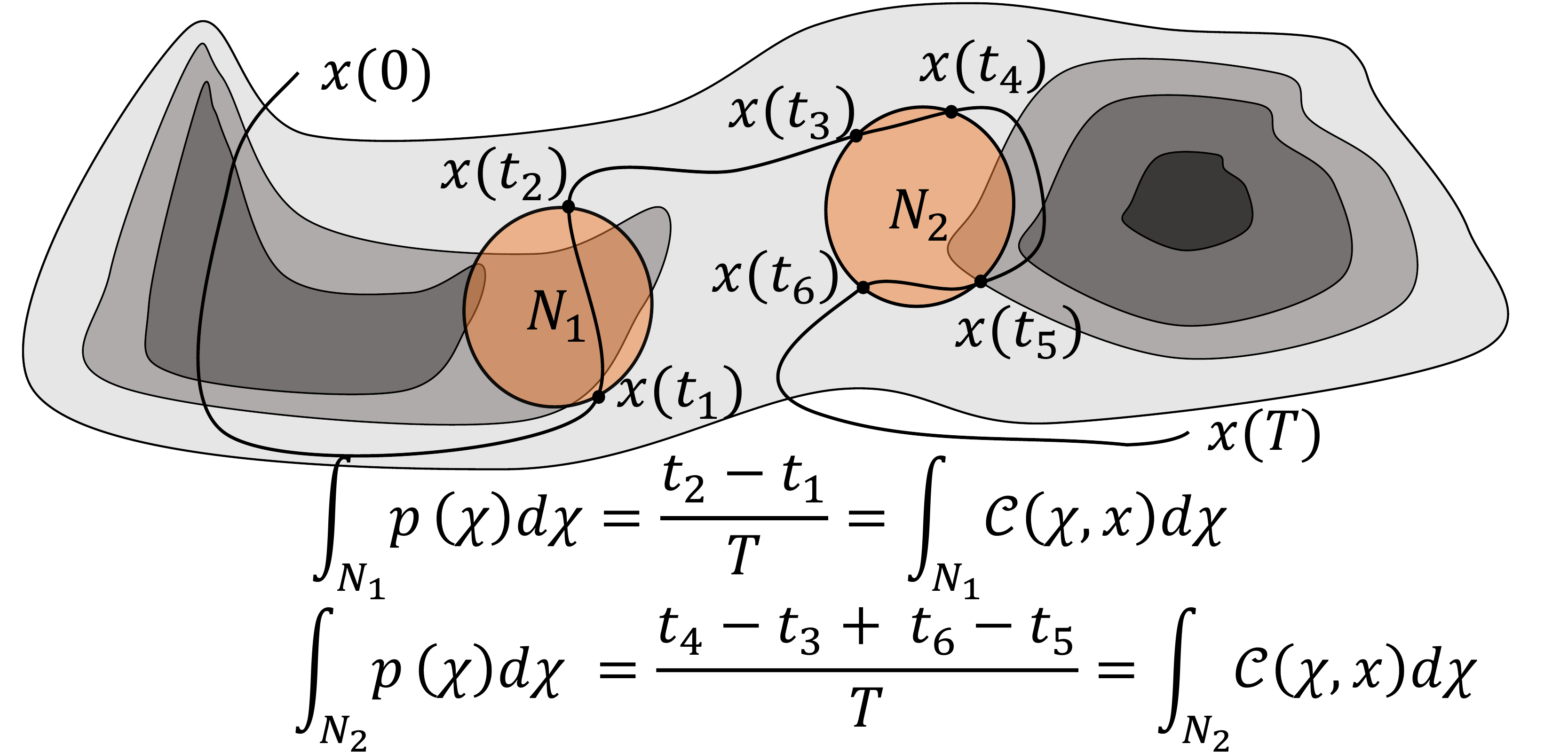}% cmpr_Bode_PRBS
\caption{Visualization of the conditions for ergodicity. }
\label{fig:erg_concept}
\end{figure}

%This section provides
A cursory overview of the ergodic trajectory planning from \cite{Mathew} is presented.
% In simple words, a system is considered as ergodic if the fraction of time spent exploring an area by system's trajectory is equal to a metric corresponding to
% the density of information of this area.
Consider a rectangular domain $\mathcal{X} := \left[0, L_{1}\right] \times \ldots \times\left[0, L_n\right]$ and an associated density of information formulated by a probability density function $p(\chi)$, with $\chi \in \mathcal{X}$. For an
horizon $T$, denote by: $x(t):[0, T] \rightarrow \mathbb{R}^n$ the state; $u(t):[0, T] \rightarrow \mathbb{R}^m$ the control; and $\dot{x}(t) = f(x(t), u(t))$ the system dynamics.  The following distribution gives information on the time-averaged statistics of a trajectory in $\mathcal{X}$
\begin{equation}\label{eq:traj_distr}
	\mathcal{C}(\chi, x)=\frac{1}{T} \int_{0}^{T} \delta(\chi- x(\tau)) d \tau,
\end{equation}
where $\delta(\cdot)$ is a Dirac delta function. Fourier series representations of $p(\chi)$ and $\mathcal{C}(\chi, x)$ are obtained by making use of the basis functions %from% using basis functions
\begin{equation}\label{eq:basis_fun}
	F_{k}(\chi)=\frac{1}{h_{k}} \prod_{i=1}^{n} \cos \left(\frac{k_i \pi}{L_{i}} \chi_{i}\right),
\end{equation}
where: $k=\left(k_{1}, k_{2}, \ldots, k_{n}\right) \in \mathbb{Z}^{n}$ denotes a multi-index that belongs to the set
$\mathcal{K}=\left\{k \in \mathbb{Z}^{n}| 0 \leq k_{j} \leq K_j\right\}$ with $K_j$ referring to the highest selected cosine harmonic; and $h_{k}$
% \begin{equation}\label{eq:norm_factor}
%	=\left(\int_{0}^{L_{1}} \ldots  \int_{0}^{L_{n}}  \prod_{i=1}^{n}\cos ^{2}\left( \frac{k_{i} \pi}{L_{i}} \chi_{i}\right) \mathrm{d} \chi_{1}
%	\ldots \mathrm{d} \chi_{n}\right)^{1 / 2},
%\end{equation}
is a normalization factor that guarantees that the Fourier basis functions $F_{k}$ have unit norm. The Fourier coefficients of the spatial and time-averaged distributions, respectively $p_k$ and $c_{k}$, can then be obtained through a standard inner product over the exploration domain $\mathcal{X}$. %as follows
% \begin{equation}
%	p_{k}=\int_{\mathcal{X}} p(\chi) F_{k}(\chi) d \chi,
%\end{equation}
%while the trajectory Fourier coefficients $c_k$ are computed similarly as
% \begin{equation}\label{eq:traj_coef}
%	c_{k} (x)=\int_{\mathcal{X}} \mathcal{C}(\chi, x) F_{k}(\chi) d \chi = \frac{1}{T} \int_{0}^{T} F_{k}(x(\tau)) d \tau .
%\end{equation}
Figure \ref{fig:erg_concept} shows the conditions for ergodicity of a trajectory $x$ with respect to subsets $N_1$ and $N_2$ (level sets of $p(\chi)$ are also reported).

In the multi-agent case, superscripts are used to denote the agents index (e.g. $x^{(j)}$ is the state of the $j$-th agent) and
%$\dot{x}^{(j)} = f(x^{(j)}, u^{(j)})$ the dynamics .
$x = \left( x^{(1)}, \ldots, x^{(N)} \right)$
and $u = \left( u^{(1)}, \ldots, u^{(N)} \right)$ denote the stacked state and control vectors, where $N$ is the number of agents. A global trajectory Fourier coefficient
$C_k (x)$ that accounts for the effect of all agents trajectories is defined as
\begin{equation}
	C_{k}(x) =\frac{1}{N} \sum_{j=1}^{N} c_{k}\left(x^{(j)}\right).
\end{equation}
%To quantify the difference between the spatial distribution and the time-averaged trajectory distribution, we define a
The \emph{shared} ergodic metric, playing a key role in the problem formulation, captures the ergodicity of the multi-agent configuration as the weighted squared
difference between the spatial and the time-averaged trajectory distributions via the respective Fourier coefficients
%distance of   $p_k$ and  $C_k$
 \begin{equation}\label{eq:ergodic_metric}
	\mathcal{E}=\sum_{k \in \mathcal{K}} \Lambda_{k}\left(C_{k}(x) -p_{k}\right)^{2},
\end{equation}
where $\Lambda_{k}=\left(1+\|k\|^{2}\right)^{-\frac{\lambda+1}{2}}$ and $\lambda \leq n$ is the number of exploratory variables in $\mathcal{X}$.

%DECENTRALIZED ERGODIC EXPLORATION
\section{DECENTRALIZED ERGODIC EXPLORATION}\label{sIII}
\Resp{This section presents the proposed decentralized algorithm for ergodic trajectory planning, which extends the single-agent problem formulation from \cite{Miller} to a multi-agent setting described by undirected and connected network topologies.}

%Also, we
Let us denote by $\xi^{(j)} :=(\alpha^{(j)}, \mu^{(j)})$ a \emph{planning trajectory} of the $j$-th agent; that is, a state-control pair (with $\alpha(t):[0, T] \rightarrow \mathbb{R}^n$ and $\mu(t):[0, T] \rightarrow \mathbb{R}^m$) which do not necessarily satisfy the system dynamics.
Let us also denote by $\eta^{(j)} := (x^{(j)}, u^{(j)})$ a feasible trajectory of the $j$-th agent, that is a state-control pair belonging to the manifold $\mathcal{T}$ of trajectories satisfying $\dot{x}^{(j)} = f(x^{(j)}, u^{(j)})$. See \cite{Hauser_manifold} for a formal characterization of $\mathcal{T}$. We also define $\xi = \left( \xi^{(1)}, \ldots, \xi^{(N)} \right)$ and $\eta = \left( \eta^{(1)}, \ldots, \eta^{(N)} \right)$ the augmented vectors consisting of the trajectories of all agents.

The problem is formulated as the minimization of an objective function featuring three contributions
%an objective function that consists of 1) an ergodic metric term, 2) a control energy index and 3) an inter-agent distance cost term as follows
 \begin{equation}\label{eq:dec_obj}
	\begin{aligned}
	&J(\xi) =  \underbrace{q \sum_{k \in \mathcal{K}}^{K} \Lambda_{k}\left(C_{k}(\alpha)-p_{k}\right)^{2}}_\text{ergodicity}
	+ \underbrace{\int_{0}^{T}  \sum_{j=1}^{N} \frac{1}{2} \left\Vert\mu^{(j)}(\tau)\right\Vert^{2}_{R(\tau)}d \tau }_\text{control energy} \\
	&+ \underbrace{\int_{0}^{T} \sum_{j=1}^{N} \sum_{\ell=j+1}^{N} \frac{1}{r_{j \ell}(\tau)+\frac{1}{2}\left\Vert \alpha^{(j)}(\tau)-\alpha^{(\ell)}(\tau) \right\Vert^{2}_{W_{j\ell}(\tau)}} d \tau}_\text{inter-agent distance},
	\end{aligned}
\end{equation}
with the following design parameters: $q>0 \in \mathbb{R}$ penalizes ergodicity; $R(\cdot) \succeq 0 \in \mathbb{S}^{m \times m}$ is a (time-varying) penalty for the control energy; $r_{j\ell}(\cdot)>0 \in \mathbb{R}$ is a (time-varying) penalty for the inter-agent distance; $W_{j\ell}(\cdot) \succeq 0 \in \mathbb{S}^{n\times n}$ is a linear transformation that allows distance between two agents positions to be computed.
\Resp{The inter-agent distance cost
encourages each agent, via the choice of $r_{j\ell}$, to perform exploration while avoiding collisions with the others.
%In particular, the inter-agent distance integrated in $[0, T]$ is penalized.
}
The convention $||x||\Resp{^2}_Q=x^\top Q x$ will be used throughout.

The goal is to determine planning trajectories $\xi ^{(j)} \in \mathcal{T}$ that minimize (\ref{eq:dec_obj}).
%\begin{equation}
%	\xi^* = \underset{\xi^{(j)} \in \mathcal{T}}{\arg \min } J \left( \xi^{(1)} , \ldots, \xi^{(N)}\right) .
%\end{equation}
Following \cite{Hauser}, the nonlinear constraint imposed by the trajectory manifold $\mathcal{T}$ is removed by making use of a projection operator $\mathcal{P}$, which
maps planning trajectories $\xi$ to feasible trajectories $\eta$. That is, the following optimization problem is considered
\begin{equation}\label{eq:proj_opt}
%	\xi^* = \underset{\xi^{(j)}}{\arg \min } J \left( \mathcal{P}\left( \xi^{(1)}\right), \ldots, \mathcal{P}\left(\xi^{(N)}\right) \right).
	\underset{\xi}{\min }\; J \left( \mathcal{P}\left( \xi^{(1)}\right), \ldots, \mathcal{P}\left(\xi^{(N)}\right) \right),
\end{equation}
and the optimal trajectory is taken as $\eta^*=\mathcal{P}(\xi^*)$, where $\xi^*$ is a minimizer of (\ref{eq:proj_opt}). Problem (\ref{eq:proj_opt}) is non-convex, and thus a local minimizer $\xi^*$ is sought. To this end, given information density Fourier coefficients $p_k$ and initial candidate trajectories $\xi_0$, (\ref{eq:proj_opt}) is solved via an iterative steepest descent algorithm, whereby each agent optimizes, in parallel and only using information from neighbouring agents, its own trajectory with the goal of minimizing the global cost (\ref{eq:dec_obj}). The iterative algorithm consists of four steps, which are the topic of the next subsections.
%\begin{enumerate}
% 	\item Each $j$-th agent seeks a descent direction $\zeta_i ^{(j)}$ to minimize locally the objective $J$ using estimations for trajectories of other agents
% 		given by its neighbourhood.
% 	\item Line search is used to determine step size $\gamma_i ^{(j)}$  to calculate the update $\xi^{(j)} + \gamma_i^{(j)} \zeta_i^{(j)}$ towards
% 		 $\zeta_i^{(j)}$ direction.
% 	\item Operator $\mathcal{P}(\cdot)$ projects the previous update into the trajectory manifold $\mathcal{T}$.
%	\item Each $j$-th agent employs the neighbourhood topology to update estimations for the rest of the agents .
% \end{enumerate}

%Steepest Descent
\subsection{Step one: steepest descent}
All agents maintain estimates of the trajectories of all other agents, thus denote by $\left\{\eta^{(\ell)}\right\}_{j}$ the feasible trajectory of the $\ell$-th agent estimated by the $j$-th agent, and by $\left\{\eta \right\}_{j} = \left( \left\{\eta^{(1)}\right\}_{j}, \ldots, \left\{\eta^{(N)}\right\}_{j}  \right)$ the stacked vector. Same convention is used for the planning trajectories $\xi$.

At iteration $i$, each agent determines its descent direction $\zeta_i ^{(j)}$ to minimize locally the objective function $J$.
To this end, the objective function $J$ is approximated by a second-order Taylor expansion around  $\xi_i ^{(j)}$ %in a Quasi-Newton philosophy, as follows
\begin{equation}\label{eq: opt_decentr_alg}
	\underset{\zeta^{(j)}_i \in T_{\xi^{(j)}_{i}} \mathcal{T}}{\min }D_{ \xi_{i}^{(j)}}  J\left(\left\{\xi_{i}\right\}_j\right) \circ \zeta_i ^{(j)}+\frac{1}{2}\left(\zeta_i ^{(j)}, \zeta_i ^{(j)}\right),
\end{equation}
where $T_{\xi^{(j)}_{i}} \mathcal{T}$ is the tangent trajectory manifold of $\mathcal{T}$, $D_{ \xi_{i}^{(j)}}  J\left(\left\{\xi_{i}\right\}_j\right) \circ
\zeta_i ^{(j)}$ is the first Frechet directional derivative and $\left(\zeta_i ^{(j)}, \zeta_i ^{(j)}\right)$ is a quadratic approximation of the second Frechet
directional derivative. The constraints force the direction $ \zeta^{(j)}_i$ to lie on the tangent space of the trajectory manifold $\mathcal{T}$ \cite{Hauser}.
%The following analysis treats each term of (\ref{eq: opt_decentr_alg}) explicitly to reveal the underlying structure of the optimization problem.
The direction $\zeta_i ^{(j)} = \left(z_i ^{(j)}, v_i ^{(j)} \right)$ is divided into state direction $z_i ^{(j)}$ and control direction $v_i ^{(j)}$.
For clarity, we omit the iteration subscript $i$ and assume that a feasible trajectory $\eta ^{(j)} $ is available (this is natural considering that projection operator $\mathcal{P}$ is applied at the end of each iteration).

\subsubsection{First Frechet directional derivative }
The first Frechet directional derivative $D_{ \eta^{(j)}}  J\left(\left\{\eta \right\}_j\right) \circ \zeta^{(j)}$can be written as
\begin{equation}\label{eq:Frechet}
 \begin{aligned}
%	 &D_{ \eta^{(j)}}  J\left(\left\{\eta \right\}_j\right) \circ \zeta^{(j)}   = \\
	& D_{x^{(j)}} J\left(\left\{\eta\right\}_j\right) \circ z^{(j)} + D_{u^{(j)}} J\left(\left\{\eta\right\}_j\right) \circ v^{(j)} \\
	&= \int_{0}^{T} a_j (\tau)^{\top}  z^{(j)}(\tau) d \tau + \int_{0}^{T}b_j(\tau)^{\top}  v^{(j)}(\tau) d \tau ,
\end{aligned}
\end{equation}
where
\begin{equation}
\begin{aligned}
a_j (\tau) &=  \sum_{k \in \mathcal{K}} \frac{2q\Lambda_{k}}{N\cdot T}\left(C_{k}\left(\left\{x\right\}_j\right)-p_{k}\right) \nabla_{x^{(j)}} F_{k}(x^{(j)}(\tau))\\
 &+  \sum_{i=1}^{N} \frac{-W_{j i}(\tau)\left(x^{(j)}(\tau)-\left\{x^{(i)}(\tau)\right\}_j\right)}{\left(r_{j i}(\tau)+\frac{1}{2}\left\Vert x^{(j)}(\tau)-\left\{x^{(i)}(\tau)\right\}_j\right\Vert^2_{W_{ji}(\tau)} \right)^{2}} .
\end{aligned}
\end{equation}
and
\begin{equation}
 b_j(\tau)= R(\tau)  u^{(j)}(\tau) .
 \end{equation}

\subsubsection{Second Frechet directional derivative}
A quadratic approximation of $(\zeta ^{(j)}, \zeta ^{(j)})$ can be obtained as
\begin{equation}\label{eq:Quad}
\int_{0}^{T}\left( ||z^{(j)}(\tau)||\Resp{^2}_{Q_{n}(\tau)}+ ||v^{(j)}(\tau)||\Resp{^2}_{R_{n}(\tau)}\right) d \tau + ||z^{(j)}(T)||\Resp{^2}_{P_{1 n}}
\end{equation}
where $Q_n (\cdot) \in \mathbb{R}^{n \times n}$ and $P_{1n} \in \mathbb{R}^{n \times n}$  are positive semi-definite
and $R_n (\cdot) \in \mathbb{R}^{m \times m}$ is positive definite.

\subsubsection{Optimization}
Using (\ref{eq:Frechet}) and (\ref{eq:Quad}), the optimization problem to find the descent direction for the trajectory update can be formulated as
\begin{equation}\label{eq:dir_opt_dec_alg}
	\begin{aligned}
	&\underset{ \left(z^{(j)}, v^{(j)}\right)}{ \min }  \int_{0}^{T} \left( a_j(\tau)^{\top} z^{(j)}(\tau)+b_j(\tau)^{\top} v^{(j)}(\tau) \right. \\
&\left. +\frac{1}{2}||z^{(j)}(\tau)||\Resp{^2}_{Q_{n}(\tau)}+\frac{1}{2}||v^{(j)}(\tau)||\Resp{^2}_{R_{n}(\tau)} \right) d \tau +\frac{1}{2} ||z^{(j)}(T)||\Resp{^2}_{P_{1 n}}\\
%&+\frac{1}{2} z^{(j)}(\tau)^{^{\top}} Q_{n}(\tau) z^{(j)}(\tau) +\frac{1}{2} v^{(j)}(\tau) R_{n}(\tau) v^{(j)}(\tau) \large) d \tau  \\
%	&+\frac{1}{2} z^{(j)}(T)^{\top} P_{1 n} z^{(j)}(T) \\
	&\textit{s.t.}\quad \dot{z}^{(j)}=\frac{\partial f}{\partial x^{(j)}} z^{(j)}+\frac{\partial f}{\partial u^{(j)}} v^{(j)}.
	\end{aligned}
\end{equation}
where the linearized dynamics constraint enforces that $\zeta^{(j)} \in T_{\xi^{(j)}} \mathcal{T}$. The optimization problem in (\ref{eq:dir_opt_dec_alg}) provides an optimal descent direction \cite{Hauser} and can be solved via differential Riccati equations \cite{Anderson_Moore}.

%Armijo line search
\subsection{Step two: Armijo line search}
To update the current planning trajectory $\xi^{(j)}_i $ with respect to the descent direction $\zeta^{(j)}_i$, %it is required to determine a step size $
%\gamma^{(j)}_i$ such that $\eta^{(j)}_{i+1} = \mathcal{P}\left(\xi^{(j)}_i  + \gamma^{(j)}_i  \zeta^{(j)}_i\right)$. For this purpose, it is
we adopt the common Armijo rule \cite{Armijo}. %, is applied and reveals the magnitude to move towards the steepest descent direction $\zeta^{(j)}_i$.
The optimal step-size is defined as
\begin{equation}\label{eq:dec_line_search}
\begin{array}{cl}
\underset{\gamma^{(j)}_{i} \in(0,1]}{\max}& \gamma^{(j)}_{i}\\
\text {s.t. }&  J(\ldots, \xi_{i}^{(j)}+ \gamma^{(j)}_{i}\zeta_i ^{(j)}, \ldots ) - J\left(\left\{\xi_{i}\right\}_j\right)\\ &
\leq \rho \cdot \gamma^{(j)}_{i}  \cdot D_{ \xi_{i}^{(j)}} J\left(\left\{\xi_{i}\right\}_j\right) \circ \zeta^{(j)}_{i}.
\end{array}
\end{equation}
%and can be interprete one still providing a decrease condition ensureslocal improvement of the updated objective function, as follows
where the user-defined parameter $\rho$ sets the desired magnitude decrease required to achieve an adequate objective improvement.
% For numerical purposes, the step-size can be parameterized as $\gamma^{(j)}_i = \beta ^h$ with $\beta \in (0, 1]$ such that the search can be done on $h \in \mathbb{N}$ \cite{kelley_iterative}.
\Resp{It is important to observe that the Armijo rule needs global information to guarantee the cost function decrease. Here, in line with the approximation proposed for step one, each agent implements a local version of the rule using the current iterates of the neighbours and the estimates of the trajectories of all other agents.
This is a heuristic, whose suboptimality we intend to analyse in future research. A starting point could be extending the local Armijo rule for strictly convex and separable cost functions with linear constraints developed in \cite{Zargham_2012_ACC} to the non-convex setting analyzed in this paper.
%This is an heuristic choice, whereas e.g. in \cite{Zargham_2012_ACC} a local Armijo update rule with guaranteed convergence was derived for strictly convex and separable cost functions with linear constraints. It is a topic for future research to quantify the suboptimality of the proposed decentralization of the update rule in the non-convex setting analyzed in this paper.
%Alternatively, the step-size can be set a-priori to a fixed small value
%This is an heuristic choice, and it is a topic for future research to quantify the suboptimality of such decentralization of the Armijo rule in the non-convex setting analyzed in this paper. In \cite{Zargham_2012_ACC} a local update rule was derived for strictly convex and separable cost functions with linear constraints.
}
\begin{figure*}[t]
	\centering
	\framebox{{ \includegraphics[scale=0.25]{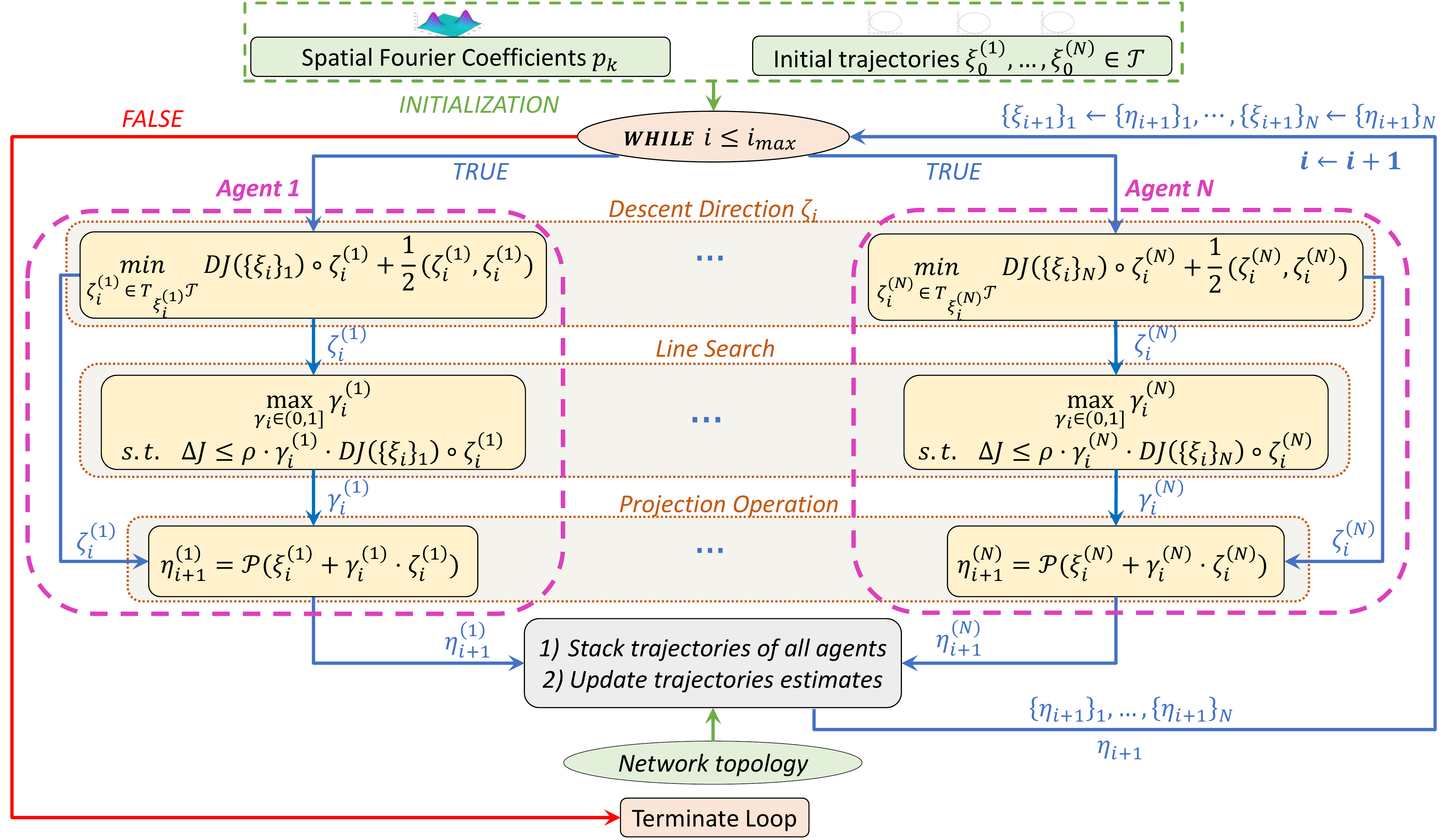}}}
	\caption{Flowchart of the proposed multi-agent decentralized ergodic trajectory optimization algorithm. }
	\label{fig:decentralized_algorithm}
 \end{figure*}
%Projection Operation
\subsection{Step three: projection on feasible trajectories manifold}

%The projection operator $\mathcal{P}$ maps any  trajectory into the trajectory manifold $\mathcal{T}$ that satisfies state space dynamics $\dot{x}^{(j)}(t) = f(x^{(j)}(t), u^{(j)}(t))$ or equivalently $\mathcal{P}: \xi =(\alpha^{(j)} , \mu^{(j)}) \mapsto \eta =(x^{(j)} , u^{(j)})$. In particular, projection operator $
%\mathcal{P}$ obtains the structure of a feedback  stabilizing compensator and turns out to be a trajectory tracking controller with the following mathematical interpretation

The projection operator $\mathcal{P}$ maps a planning trajectory into the closest one belonging to the trajectory manifold $\mathcal{T}$ made of trajectories that satisfy $\dot{x}^{(j)} = f(x^{(j)}, u^{(j)})$, that is $\mathcal{P}: \xi^{(j)} =(\alpha^{(j)} , \mu^{(j)}) \mapsto \eta^{(j)} =(x^{(j)} , u^{(j)})$. It was shown in \cite{Hauser} that the projection operator $\mathcal{P}$ can be interpreted as the following trajectory tracking controller %with the following feedback interpretation
  \begin{equation}\label{eq:proj}
	\eta^{(j)} = \mathcal{P}(\xi^{(j)}):
	\left\{\begin{array}{l} u^{(j)}=\mu^{(j)}+K^{(j)}(\alpha^{(j)}-x^{(j)}) \\
	\dot{x}^{(j)}=f(x^{(j)}, u^{(j)}) .
	\end{array}\right.
\end{equation}
The optimal controller gain $K^{(j)}$ can be computed as the solution of a finite horizon Linear Quadratic Regulator problem \cite{Anderson_Moore} applied to a linearization of the nonlinear dynamics around the trajectory $\xi^{(j)} =(\alpha^{(j)} , \mu^{(j)})$. A useful criterion for choosing the weighting matrices $R_{LQR} \in \mathbb{R}^{m\times m}$ and  $Q_{LQR} \in \mathbb{R}^{n\times n}$
is to tune them so that the projected trajectory $\eta^{(j)}$ is close to $\xi^{(j)}$ and thus the linearization gives sufficiently accurate results.
%(see more in Section \ref{sIV}).
%The LQR weighting matrices have been tuned such that changes of the trajectories $\xi_i^{(j)}$ ensure the linearization required for the projection.

%Network connectivity and estimation
\subsection{Step four: agents trajectories estimation}

Agent $j$ has only access to the trajectories optimized by the agents
%restricted information about the whole system. To encode this constraint, we introduce the notion of the
in its neighbourhood $\mathcal{N}(j)$. % for
%each $j$-th agent, that is the set of agents which interchange information directly with the $\ell$-th agent.
At the end of round $i$, the vector of feasible trajectories estimated by the $j$-th agent $\left\{\eta^{(\ell)}\right\}_{j}$ is obtained by the communication protocol below
\begin{equation}\label{eq:neighbourhood}
\left\{\eta_{i}^{(\ell)}\right\}_{j}=\left\{\begin{aligned}
	&\eta_{i}^{(\ell)} \quad  \text{if}\quad \ell \in\left\{\mathcal{N}(j)\cup\{j\}\right\} \\
	&\frac{1}{|\mathcal{N}(j)|+1}\left(\sum_{k \in \mathcal{N}(j)\cup\{j\}} \left\{\eta_{i-1}^{(\ell)}\right\}_{k}\right),\text {otherwise}.
\end{aligned}\right.
\end{equation}
%In case an agent belongs to the neighbourhood of the $j$-th agent, then the $j$-th agent has full knowledge of the current real trajectory information of this neighbour agent.
If an agent $\ell$ is not included in the neighbourhood of the $j$-th agent, then the $j$-th agent averages its estimate of agent $\ell$'s trajectories with its neighbours' estimates of agent $\ell$'s trajectory.
Figure \ref{fig:decentralized_algorithm} illustrates a flowchart of the iterative optimization algorithm.

% \begin{color}{red}
\Resp{The termination criterion consists of stopping the algorithm after a maximum number of iterations $i_{max}$ has been performed. Because the algorithm is run offline, the value of $i_{max}$ can be chosen large enough so that the \emph{ergodic reduction} metric $\mathcal{E}_r:=100\frac{\mathcal{E}_0 - \mathcal{E}_f }{\mathcal{E}_0}$, where $\mathcal{E}_f$ and $\mathcal{E}_0$ are the ergodicity at the $i_{max}$-th and initial trajectory, respectively, has a satisfactory value.
Alternatively, each agent can broadcast to its neighbours a flag when its local problem has reached the termination criterion (which can include e.g. a local estimate of $\mathcal{E}_r$ and the local directional derivative), which in subsequent iterations is rebroadcast to their (the neighbours') neighbours, and so on.
%Once an agent receives such a flag from its neighbour (say, agent $i$), it broadcasts to its own neighbours that agent $i$ has reached the termination criterion.
Once an agent has collected termination flags from all other agents, it ceases running the algorithm.
This happens $G$ iterations after all agents have reached their respective local termination criteria, where $G$ is the length of the longest path in the graph (or \emph{girth}).
%(otherwise known as the \emph{girth}).
}

%Two alternative termination criteria are employed. The first, already proposed by \cite{Miller}, uses a tolerance $\epsilon_1$ on the directional derivative. To avoid slow convergence typically associated with this criterion, and to emphasize the importance of the ergodic metric, an additional termination criteria, called \emph{ergodic reduction} $\mathcal{E}_r:=100\frac{\mathcal{E}_0 - \mathcal{E}_i }{\mathcal{E}_0}$, is introduced, where $\mathcal{E}_i$ and $\mathcal{E}_0$ are the ergodicity at the $i$-th iteration and that of the initial trajectory, respectively, with a tolerance $\epsilon_2$. %Next, we analyze in detail steps 2-5.
%\end{color}

%IMPLEMENTATION
\section{TEST-CASE PRESENTATION}\label{sIV}
The complete model description, the algorithm's parameters and the performance metrics used to evaluate the proposed algorithm are provided in this section.
 %Agent dynamics
 \subsection{Agents dynamics}
We consider the nonlinear dynamic model for the motion of the single agent used in \cite{Miller} and a time horizon $T= 3.5 \text{ sec}$.
That is
\begin{equation}\label{state_dyn}
\dot{x}^{(j)} = f(x^{(j)}, u^{(j)})=\left[\begin{array}{cc}
\cos \left(\theta^{(j)}\right) & 0 \\
\sin \left(\theta^{(j)}\right) & 0 \\
0 & 1
\end{array}\right] u^{(j)}.
\end{equation}
The state vector is $x^{(j)} = \left[ X^{(j)}, Y^{(j)}, \theta^{(j)} \right] ^{\top}$ where
$X^{(j)}$ and $Y^{(j)}$ correspond to Cartesian coordinates, while $\theta^{(j)} $ is the heading angle of the velocity vector.
%Note that while the domain $\mathcal{X}$ is by definition in the state space, only the cartesian coordinates are exploration variables. This can be easily accounted for by setting to zero the contribution of the third state to the ergodic metric.
The control input vector $u^{(j)} = \left[ \nu^{(j)}, \omega^{(j)} \right] ^{\top}$ consists of the forward velocity $\nu^{(j)}$ and the time derivative $\omega^{(j)}$ of the heading angle $\theta^{(j)}$. The initial feasible trajectories are circles with radius $\mathcal{R} = 0.05 \text{ m}$ and center randomized as discussed later.

%Exploration field
\subsection{Exploration field}
The exploration field is assumed to be a two-dimensional space $\left[0, 1\right] \times\left[0, 1\right]$  corresponding to $X$ and $Y$ Cartesian coordinates and, thus, the number of exploratory variables is $\lambda = 2$ . The information density, assigned to this field, is modelled based on a Gaussian mixture structure as follows
\begin{equation}\label{mixtures}
	p(\chi)=\sum_{i=1} w_{i} \frac{1}{\sqrt{2 \pi\left|\Sigma_{i}\right|}} e^{\left(-\frac{1}{2}\left(M\chi-\mu_{i}\right)^{\top} \Sigma_{i}^{-1}\left(M\chi-\mu_{i}\right)\right)},
\end{equation}
where: $M=\begin{bmatrix}
		1 & 0 & 0 \\
0 & 1 & 0\end{bmatrix}$ is a transformation that maps $\chi$ to the exploration variables; $\mu_i \in \mathbb{R}^2$ is the mean vector of the $i$-th mode; $\Sigma_i \in \mathbb{R}^{2\times 2}$ is the positive definite covariance of the $i$-th mode; and $w_i \in [0, 1]$ is the weight of the $i$-th mode (chosen such that $p$ integrates to 1 in the exploration domain). To investigate the ability of the algorithm to design different planning strategies as a function of the information densities, we investigate two different cases for (\ref{mixtures}).
%\begin{itemize}
%	\item
The first, named \emph{volcano}, has a dominant mode at the center of the exploration map and other minor modes peripherally.
		 The first mode has $w_1 = 0.6$, $\Sigma_1 = 0.014 \cdot \mathcal{I}_{2\times 2}$ and $\mu_{1}=\begin{bmatrix}
		0.5 & 0.5 \end{bmatrix} ^ {\top}$, whereas the rest of the modes are equally weighted with $\Sigma_i = 0.004 \cdot \mathcal{I}_{2\times 2}$,
		 $\mu_{2}=\begin{bmatrix} 0.75 & 0.5 \end{bmatrix} ^ {\top}$, $\mu_{3}=\begin{bmatrix} 0.25 & 0.5 \end{bmatrix} ^ {\top}$
                 $\mu_{4}=\begin{bmatrix} 0.5 & 0.75 \end{bmatrix} ^ {\top}$ and   $\mu_{5}=\begin{bmatrix} 0.5 & 0.25 \end{bmatrix} ^ {\top}$.
%	\item
The second, named \emph{archipelago}, has four modes with same covariance $\Sigma = 0.006 \cdot
		 \mathcal{I}_{2\times 2}$ and mean vectors $\mu_{1}=\begin{bmatrix} 0.25 & 0.25 \end{bmatrix} ^ {\top}$,
                 $\mu_{2}=\begin{bmatrix} 0.75 & 0.25 \end{bmatrix} ^ {\top}$, $\mu_{3}=\begin{bmatrix} 0.25 & 0.75 \end{bmatrix} ^ {\top}$,
                  $\mu_{4}=\begin{bmatrix} 0.75 & 0.75 \end{bmatrix} ^ {\top}$.
%\end{itemize}

The approximation of spatial and trajectory distributions is addressed through the basis functions (\ref{eq:basis_fun}), where we set $K_1 = K_2 = 10$ for $X$ and $Y$ coordinates, respectively, and $K_3 = 0$ for  $\theta$ coordinate as there is no exploration. It is noted that an increase in state dimension (\ref{state_dyn}) has no effect on the computation of the coefficients, since one would assign zero to the indexes $K_j$ associated with non-exploratory states $j$. Since the planning is done off-line, considering a higher order system would have little impact on the rest of the trajectory optimization problem. %, thus it is believed that this test-case can provide insights into the multi-agent ergodic problem.

%and the number of corresponding $K_j$
%modes. Based on that $\Lambda_k$ in (\ref{eq:basis_fun}) prioritizes low-frequency modes, $K_1 = K_2 = 10$ is  a valid option for
%$X$ and $Y$ coordinates, while we set $K_3 = 0$ as we do not explore with respect to $\theta$ coordinate.
%Tuning parameters
\subsection{Parameters and topology}

%The termination criteria we use is $\epsilon_1 = 0.001$ and we discuss $\epsilon_2$ in Section \ref{sV}.
Table \ref{tab:tuning_parameters} summarizes design parameters for the trajectory optimization algorithm.
 In the objective function, we prioritize ergodicity against control energy by adjusting accordingly the relative magnitudes of $q$ and $R$. It can also be observed that matrix $W_{jl} = W$  extracts
from the state vector elements related to Cartesian coordinates so that a Euclidean distance metric is obtained in the inter-agent cost term. The choice for $Q_n$ and $R_n$ is motivated
by a desired smooth change on the state $x$ and a more aggressive change on control $u$, respectively.
The LQR weighting matrices have been tuned according to the previously discussed criterion.
%$\begin{bmatrix} 1 & 0 & 0 \\ 0 & 1 & 0 \\ 0 & 0 & 0 \end{bmatrix}$
\begin{table}[h]
\centering
\caption{Tuning parameters}
\begin{tabular}{ccc}
\hline
                   \textbf{Type}            & \textbf{Parameters} & \textbf{Values}                                                                                             \\ \hline
                                                  	                & $q$                   & $100$                                                                                          \\ \cline{2-3}
                                                                         & $R$              & $0.03 \cdot I_{2 \times 2}$                                                            \\ \cline{2-3}
\multirow{-3}{*}{\textbf{Objective Function $J$}} & $W$                 & $\operatorname{diag}\left(1,1,0\right)$ \\ \hline
                                                                         & $Q_n$          & $450 \cdot I_{3\times 3}$                                           		           \\ \cline{2-3}
                                                                         & $R_n$          & $14.5 \cdot I_{2 \times 2}$                                                             \\ \cline{2-3}
\multirow{-3}{*}{\textbf{Quasi-Newton}}             & $P_{1n}$          & $50 \cdot I_{3 \times 3}$	   					                    \\ \hline
                                                                        & $Q_{LQR}$    & $I_{3\times 3}$                                                                              \\ \cline{2-3}
\multirow{-2}{*}{\textbf{LQR}}	                        & $R_{LQR}$     & $I_{2\times 2}$                                                                             \\ \hline
\end{tabular}
\label{tab:tuning_parameters}
\end{table}

%\subsection{Line Search Parameters}

For the line search problem (\ref{eq:dec_line_search}), we parameterize $\gamma^{(j)}_i = \beta ^h$ with $\beta = 0.99$. With this choice, a fine search on the step-size is allowed and, thus, a larger reduction of the objective function is achieved, with convergence benefits  \cite{kelley_iterative}.
Parameter $\rho$ specifies the  magnitude of reduction in the sufficient decrease condition and a typical value $\rho = 10^{-4}$ is used. Finally, \Resp{the termination criterion is $i_{max}$ $=$ $70$}.

\begin{figure}[t]
	\centering
	\framebox{{ \includegraphics[width=\columnwidth]{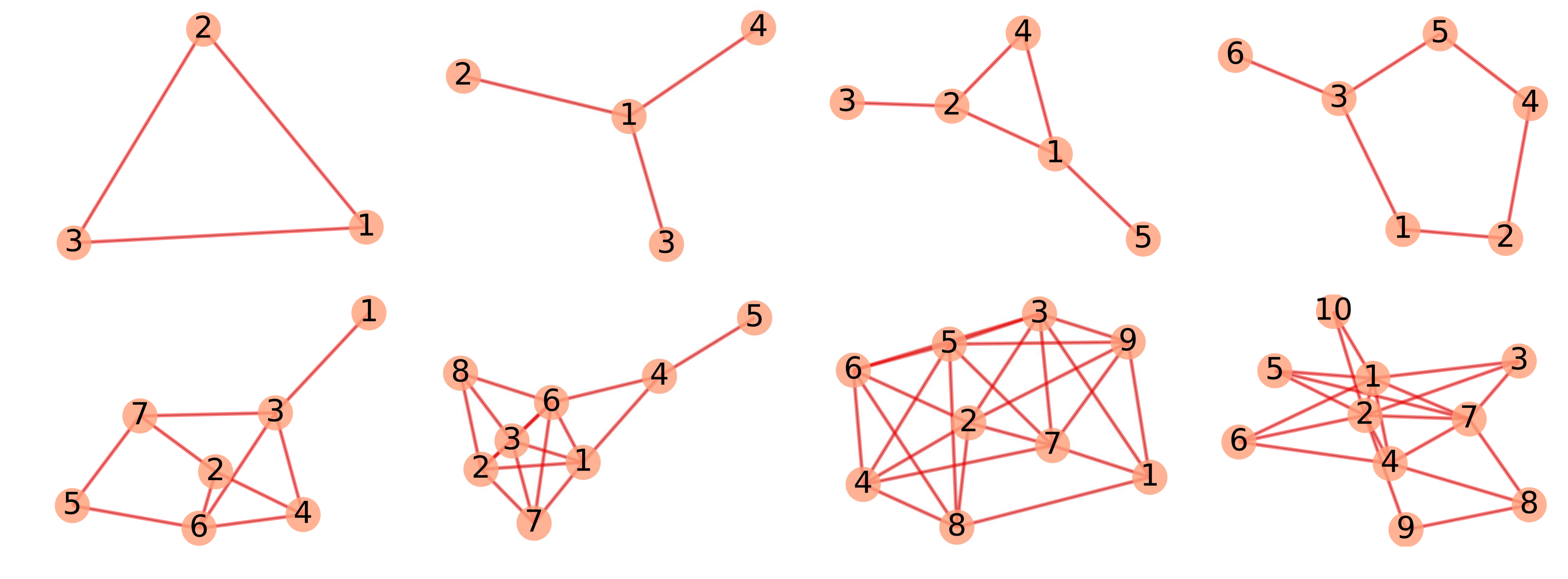}}}
	\caption{Network topology for 3$\leq N \leq$10.}
	\label{fig:erdos-renyi}
 \end{figure}

%Network Topology
%\subsection{Network Topology}
% In decentralized exploration, each agent is only aware of neighbourhood information. Hence, it is essential to define a network capable of describing the
% interconnection between the agents.
%Given the number of agents,
The analyzed network topologies for 3$\leq N \leq$10 were randomly generated and are shown in Fig. \ref{fig:erdos-renyi}.
\Resp{This topology defines the fixed neighbourhood $\mathcal{N}(\cdot)$ used in (\ref{eq:neighbourhood}) and establishing the communication constraints among the agents.}
%constructed based on Erdos–Renyi random graphs \cite{erdos}, whereby two agents are connected to each other randomly with a probability $0.5$.  This procedure for generating graphs can be used to assess stochastic performance of the proposed algorithm with respect to different network topologies. In the analyses shown in the next section, we instead always consider the network topologies depicted in Fig. \ref{fig:erdos-renyi}, and analyze stochasticity with respect to randomly sampled initial trajectories.

\subsection{Performance metrics}
%To evaluate t
The performance of the algorithm is investigated based on the four metrics defined below (the first two are global, while the last two are at agent-level):
\begin{itemize}
	\item Optimal temporal ergodic metric \cite{Mathew}:
		 \begin{equation}\label{eq:opt_ergodic_metric}
		\mathcal{E}_{opt}(t)=\sum_{k \in \mathcal{K}} \Lambda_{k}\left(C_{k}\left(x^{*}, t\right) -p_{k}\right)^{2},
		\end{equation}
	 where $C_{k}(x, t)=\frac{1}{N} \sum_{j=1}^{N} \frac{1}{t} \int_{0}^{t} F_{k}\left(x^{(j)}(\tau)\right) d \tau$. % $\mathcal{E}_{opt}(t)$
	 This is a time-dependent feature as it evaluates ergodicity with respect to any time instant in the horizon.

	\item Completion task time \cite{voronoi,nestmeyer2017decentralized}: it is generally used to indicate the time $t_{CTT}$ at which exploration can terminate without loss of beneficial information. In this work, we define it as the first time instance $t \in [0, T]$ when
$\frac{\mathcal{E}_{opt}(t=0) - \mathcal{E}_{opt}(t) }{\mathcal{E}_{opt}(t=0)} \geq \epsilon_{opt}$
\Resp{, where $\epsilon_{opt}$ is defined in Section \ref{sV}.}
%in which the optimal
% 			temporal ergodicity reduction $\frac{\mathcal{E}_{opt}(t=0) - \mathcal{E}_{opt}(t) }{\mathcal{E}_{opt}(t=0)} 100$
%			 exceeds a given tolerance $\epsilon_{opt}$.

	\item Control energy per agent:
	%	\begin{equation} 		\sqrt{\int_{0}^{t_{CTT}} \left[\nu^{(j)}(\tau)^{2} + \omega^{(j)}(\tau)^{2}\right] d \tau}.		\end{equation}
	$\sqrt{\int_{0}^{t_{CTT}} u^{(j)}(\tau)^{2} d \tau}.$
%	         It reveals the energy consumption required for an agent to implement the proposed control plan up to $t_{CTT}$.

	\item  Traveled Distance:
%		\int_{0}^{t_{CTT}} \sqrt{\left(\frac{dX^{(j)}}{dt}\right)^{2}+\left(\frac{dY^{(j)}}{dt}\right)^{2}} d t.
$\int_{0}^{t_{CTT}}\lvert \nu^{(j)}(\tau) \rvert d \tau.$
%\begin{equation} \int_{0}^{t_{CTT}}\lvert \nu^{(j)}(\tau) \rvert d \tau. \end{equation}
	%	 An agent with a large traveled distance \cite{nestmeyer2017decentralized} is more sensitive to obstacle crashes
	%	 or unpredictable  failures.
	
\end{itemize}

 \begin{figure*}[t]
     \centering
%\framebox{
     \begin{subfigure}{0.31\textwidth}
         \centering
         \includegraphics[trim={0cm 0 0cm 0},clip,scale =0.3]{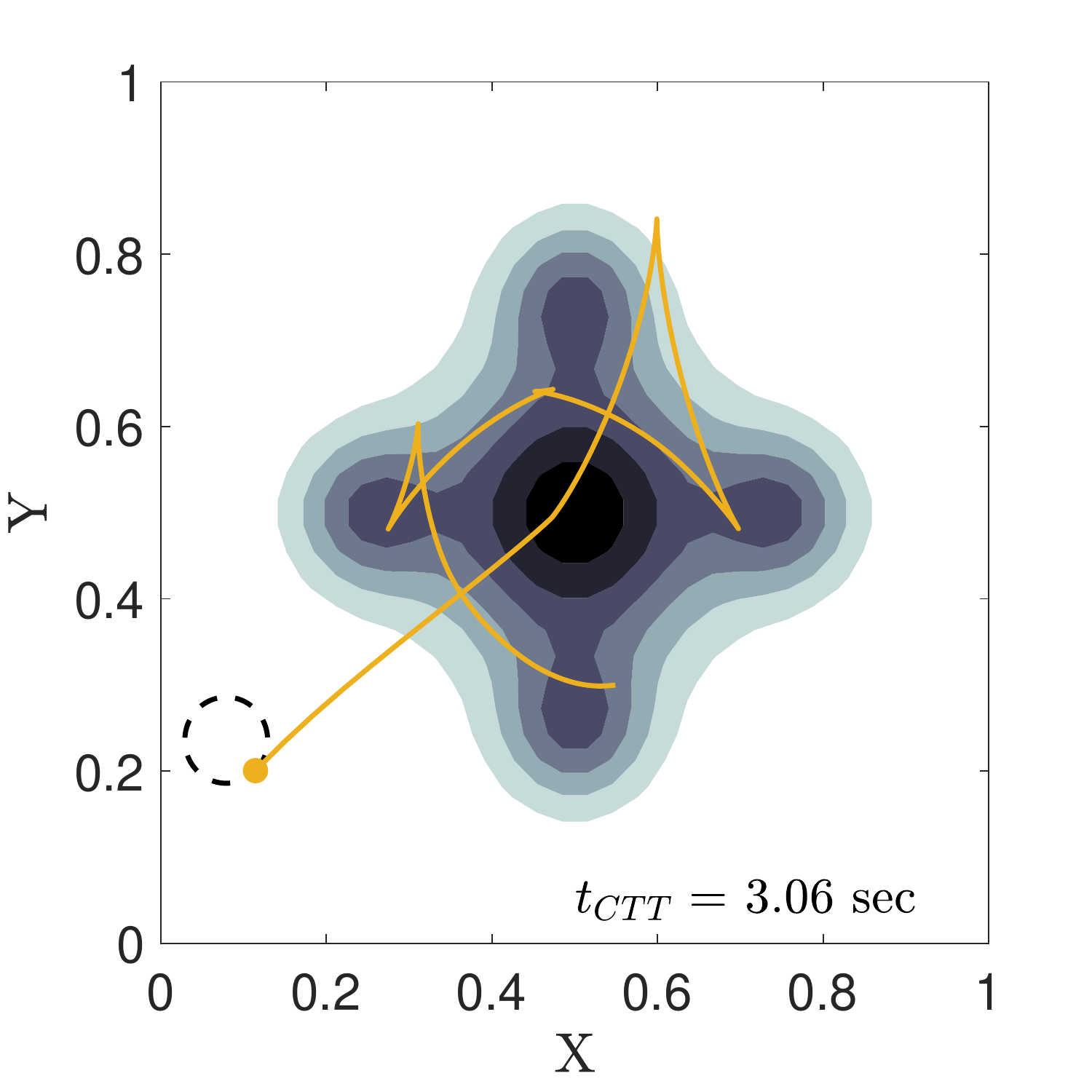}
         \caption{ }
         \label{fig:vol_1}
     \end{subfigure}
     \hfill
        \begin{subfigure}{0.31\textwidth}
         \centering
         \includegraphics[trim={0cm 0 0cm 0},clip,scale =0.3]{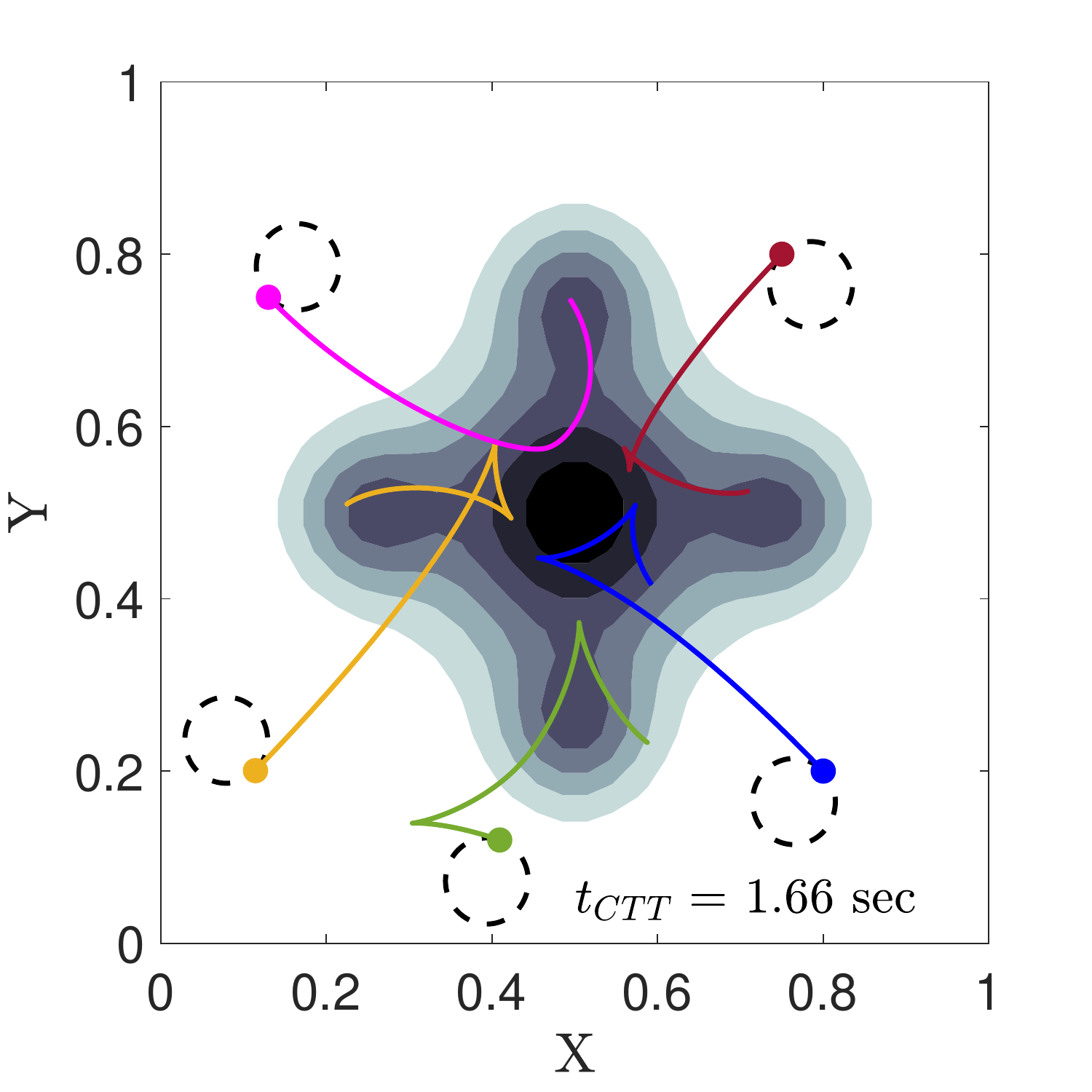}
         \caption{ }
         \label{fig:vol_2}
     \end{subfigure}
      \hfill
      \begin{subfigure}{0.35\textwidth}
         \centering
         \includegraphics[trim={0cm 0 0cm 0},clip,scale =0.45]{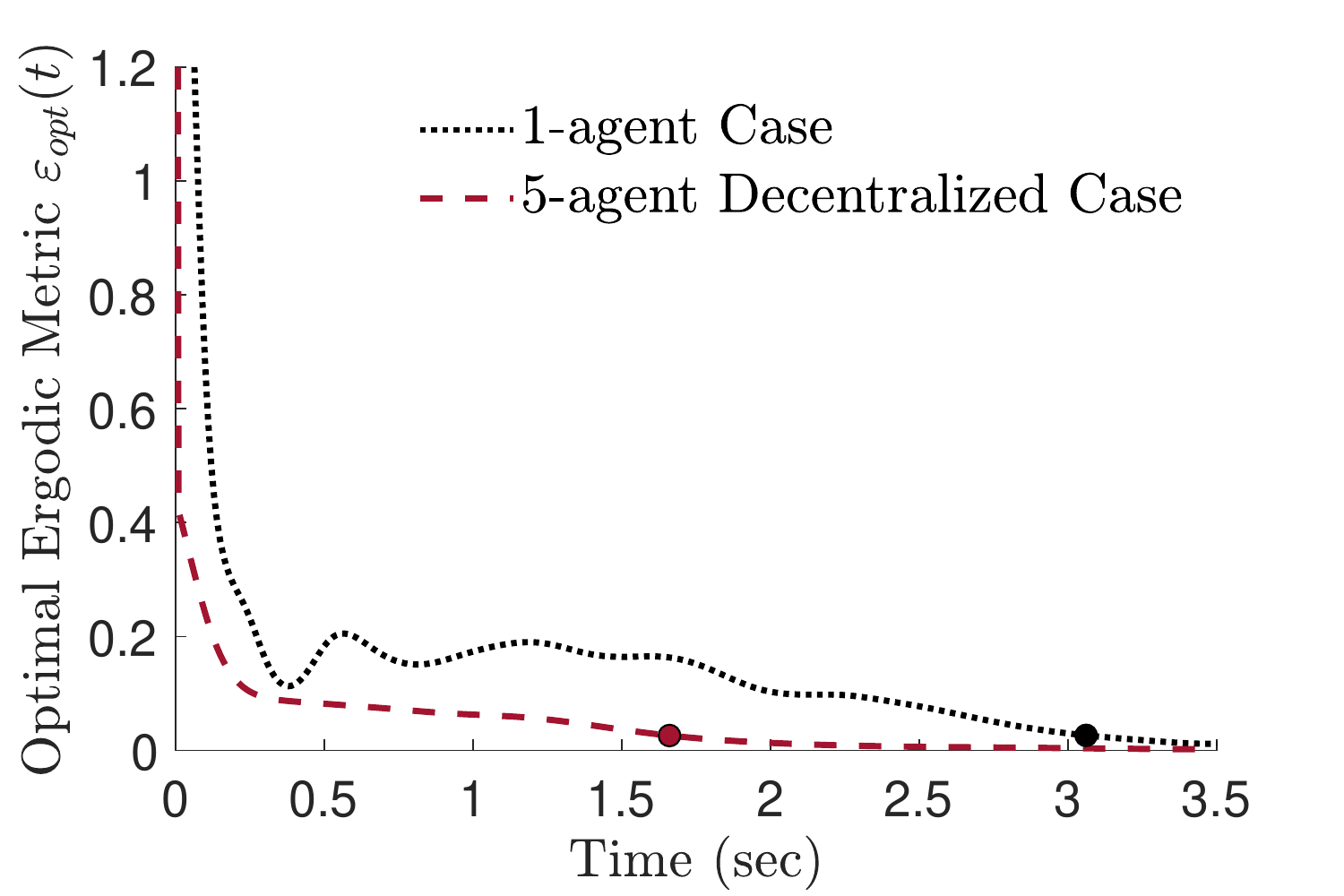}
         \caption{ }
         \label{fig:vol_3}
     \end{subfigure}
%}
     \caption{Volcano case: State trajectory $\left(X(t), Y(t)\right)$ for single agent (a) and five agents (b) from $t=0$ until $t=t_{CTT}$; (c) $\mathcal{E}_{opt}(t)$ with respect to time. Dashed lines in (a) and (b) show the initial trajectory. \Resp{Circular markers in (c) indicate $t_{CTT}$}.}
     \label{fig:vol}
\end{figure*}

 \begin{figure*}[t]
%\framebox{
     \begin{subfigure}{0.31\textwidth}
         \centering
         \includegraphics[trim={0cm 0 0cm 0},clip,scale =0.3]{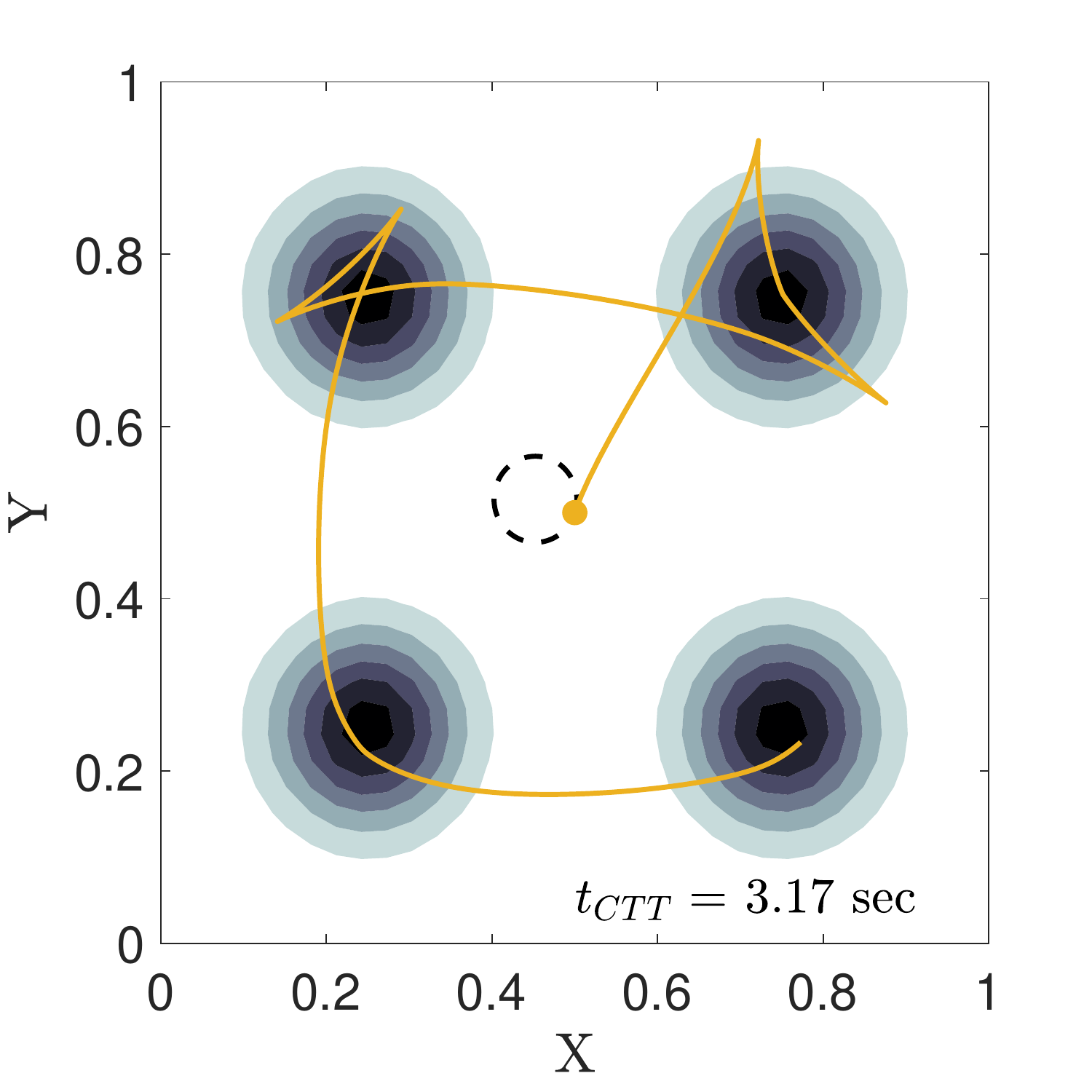}
         \caption{ }
         \label{fig:isl_1}
     \end{subfigure}
     \hfill
        \begin{subfigure}{0.31\textwidth}
         \centering
         \includegraphics[trim={0cm 0 0cm 0},clip, scale =0.3]{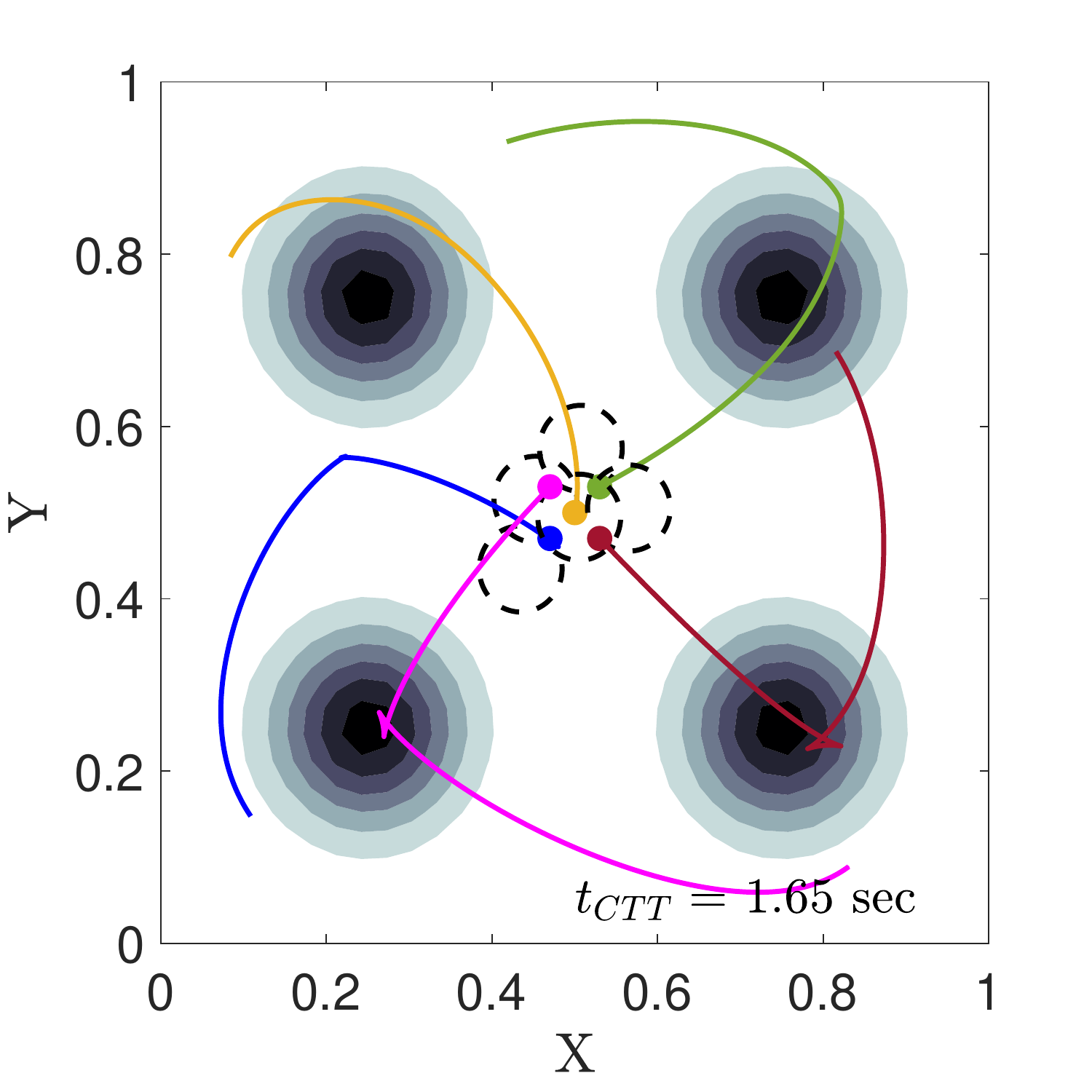}
         \caption{ }
         \label{fig:isl_2}
     \end{subfigure}
      \hfill
      \begin{subfigure}{0.35\textwidth}
         \centering
         \includegraphics[trim={0cm 0 0cm 0},clip,scale =0.45]{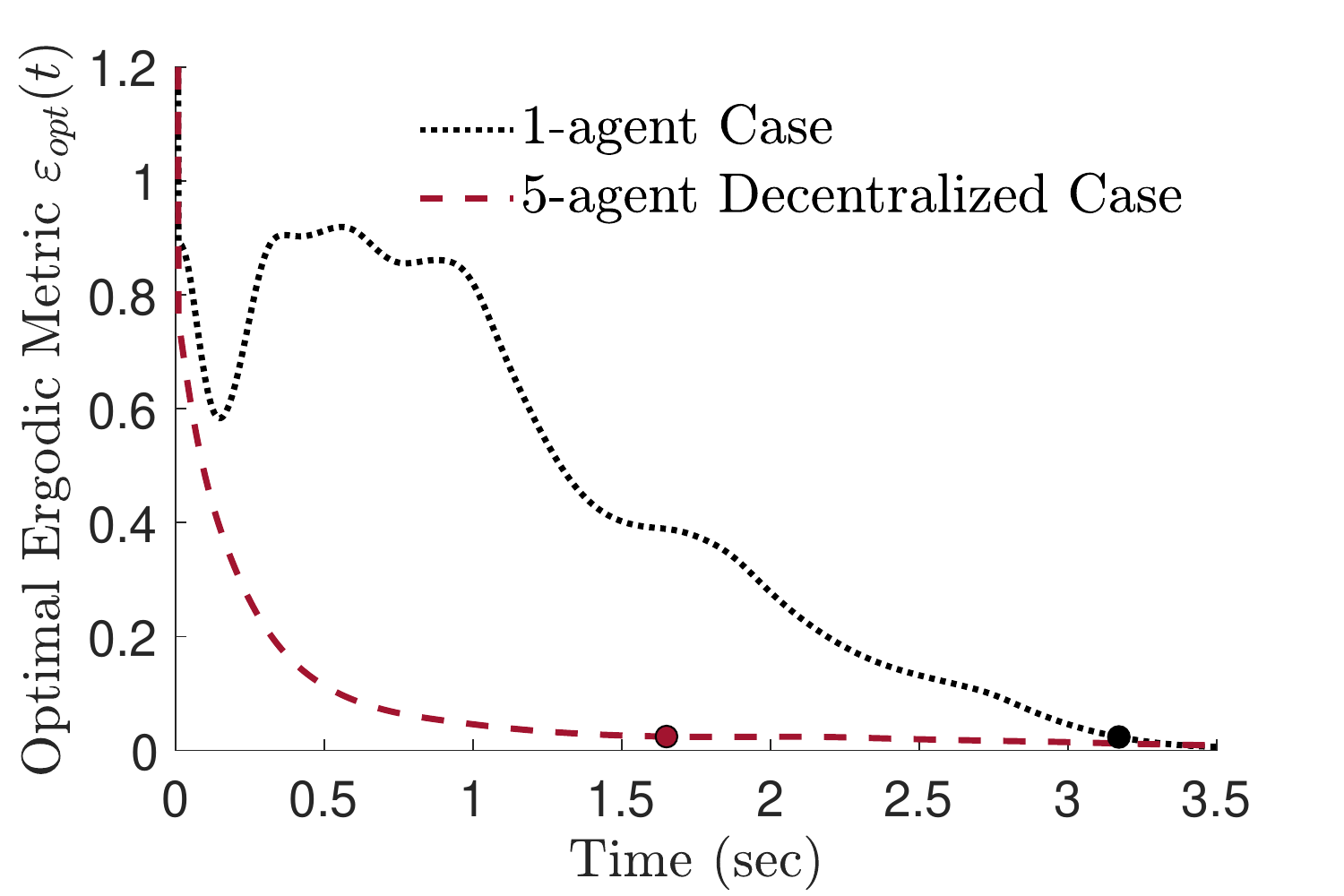}
         \caption{ }
         \label{fig:isl_3}
     \end{subfigure}
%}
     \caption{Archipelago case: State trajectory $\left(X(t), Y(t)\right)$ for single agent (a) and five agents (b) from $t=0$ until $t=t_{CTT}$; (c) $\mathcal{E}_{opt}(t)$ with respect to time. Dashed lines in (a) and (b) show initial trajectory. \Resp{Circular markers in (c) indicate $t_{CTT}$}.}
     \label{fig:isl}
\end{figure*}

 %RESULTS
 \section{RESULTS}\label{sV}
 The implementation of the algorithm used to generate the results shown in this section are available at this repository \cite{repo_ETH}.
 %Ergodic optimal trajectories
 \subsection{Optimal trajectories}
A comparison between the optimal trajectories in the Cartesian coordinates for single-agent ($a$) and 5-agent cases ($b$) with respect to the two information density distributions discussed earlier is presented in Figures \ref{fig:vol} and \ref{fig:isl}. The figures also show the optimal ergodic metric $\mathcal{E}_{opt}$ ($t$) as a function of the time inside the horizon $t$. The inter-agent distance penalty term and  optimal temporal ergodicity reduction tolerance are  set to $r = 1$ and $\epsilon_{opt} = 99.5$, respectively. The initial trajectories used in the optimization are represented in dashed line, while the system's state initial condition $x_0$ at $t=0$ is picked in low interest regions. % in order to enable unbiased comparisons between different scenarios.

Figure \ref{fig:vol_1} shows how the single-agent attempts to cover the whole region by spending more time in regions of higher interest, according to the ergodic principle. However, a high amount of time is required to cover efficiently the domain as shown by the completion task time value $t_{CTT} = 3.06 \text{ sec}$. On the other hand, in the multi-agent case in Fig. \ref{fig:vol_2}, the agents collaborate and split the field into exploration sub-domains.  This is enabled by the use of a shared ergodic metric $\mathcal{E}$ in the objective function (\ref{eq:dec_obj}) that frames exploration
as a common group task among the agents.
%Each agent explores its own sub-domain independently to the rest of the agents.
As a result, a more time-efficient exploratory mission compared to the single-agent case is accomplished (note that $t_{CTT} = 1.66 \text{ sec}$ value in Fig. \ref{fig:vol_2}).
% To provide more insights to the comparison between single-agent and multi-agent case, we present in
Figure \ref{fig:vol_3}, which displays the time-dependent optimal ergodic metric $\mathcal{E}_{opt}$, shows another distinctive feature of the multi-agent solution.
 %during the final iteration for  Volcano distribution.
Namely, $\mathcal{E}_{opt}$ is monotonic with respect to time in the 5-agent case, unlike in the single-agent case.
%It is clearly seen that the single-agent acquires initially an undershoot that degrades performance and raises $\mathcal{E}_{opt}(t)$ and, then obtains  a decreasing trend up to the time horizon $T$.
This can be explained by observing that the single-agent in Fig. \ref{fig:vol_1} displays an initially efficient area coverage around high-valued information regions, but afterwards
 it spends a fraction of time in a low-valued information. This is necessary to move towards different high-valued regions, and determines the temporary increase of $\mathcal{E}_{opt}$.
  %region shaping a  figure-triangular-like  curve.
On the contrary, in the multi-agent case this inefficient part can be avoided by leveraging the possibility to optimize over multiple trajectories and thus distribute exploration  to maximize the information reward. This feature, which had not been previously observed to the best of the authors knowledge, points out an additional benefit of the multi-agent configuration. Indeed, if $\mathcal{E}_{opt}$ is monotonic, one is guaranteed to improve, in an ergodic sense, on learning the exploration field as time proceeds.
%which can be a useful property in practical applications. % This can be helpful, for example, if a mission has to suddenly by done on a shorter horizon $T$ (e.g. in case of emergency) with limited repercussions. Instead, in the single-agent case it becomes crucial the in on-line applications. Instead, when $\mathcal{E}_{opt}$ is non-monotonic,
%having curve portrays a much smoother behavior that decreases drastically $\mathcal{E}_{opt}(t)$, showing a high-quality ergodic area coverage.
%The circular markers indicate the completion task time that, as observed earlier, is much lower for the multi-agent case. % compared to the single-agent.
%This demonstrates that multi-agent configuration offers an efficient quick exploration against single-agent case.
 \begin{figure*}[t]
     \centering
  %   \framebox{
     \begin{subfigure}{0.45\textwidth}
         \centering
         \includegraphics[trim={0cm 0 0cm 0},clip, scale =0.38]{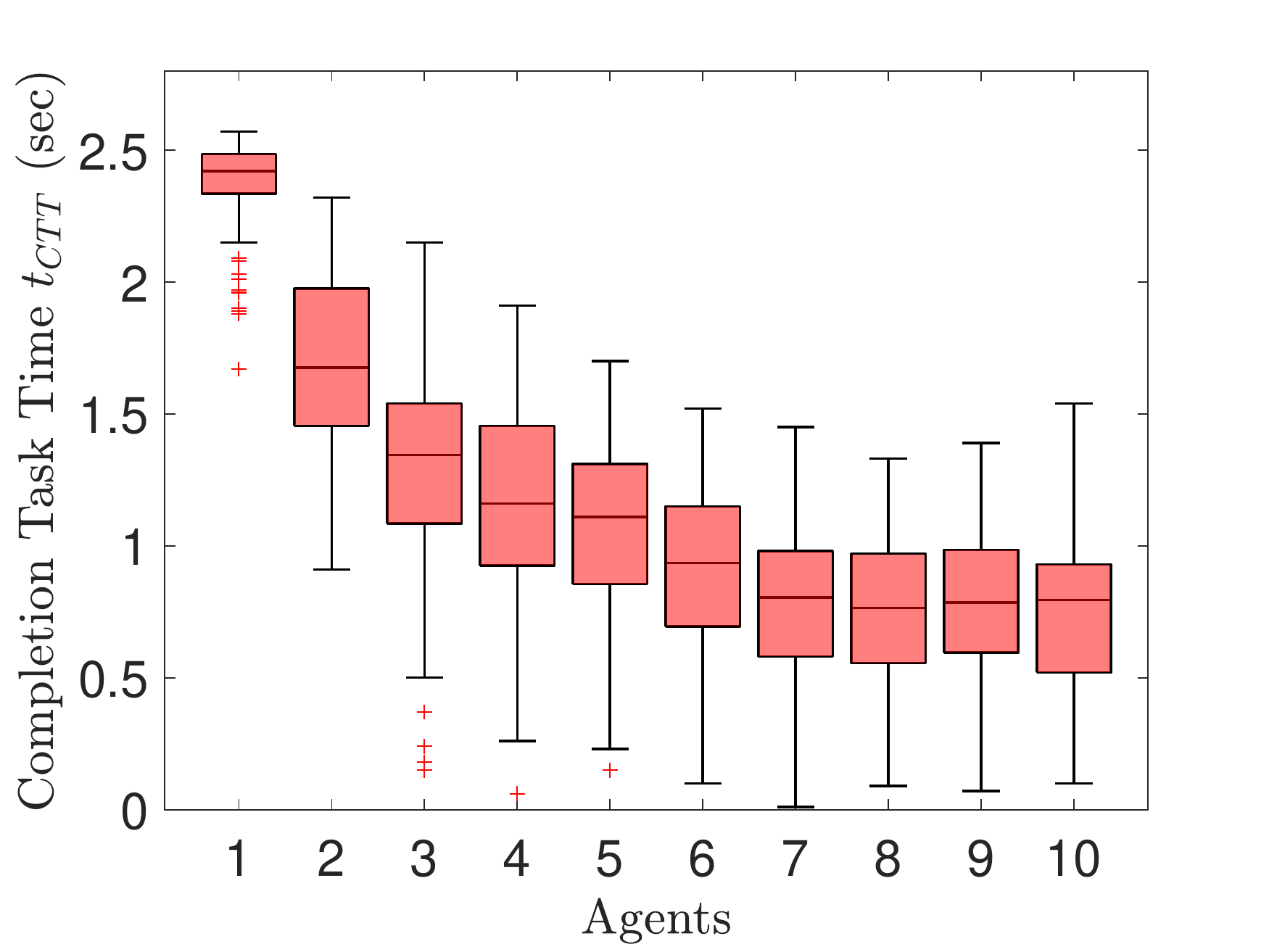}
         \caption{ }
         \label{fig:bp1a}
     \end{subfigure}
     \hfill
     \begin{subfigure}{0.45\textwidth}
         \centering
         \includegraphics[trim={0cm 0 0cm 0},clip, scale =0.38]{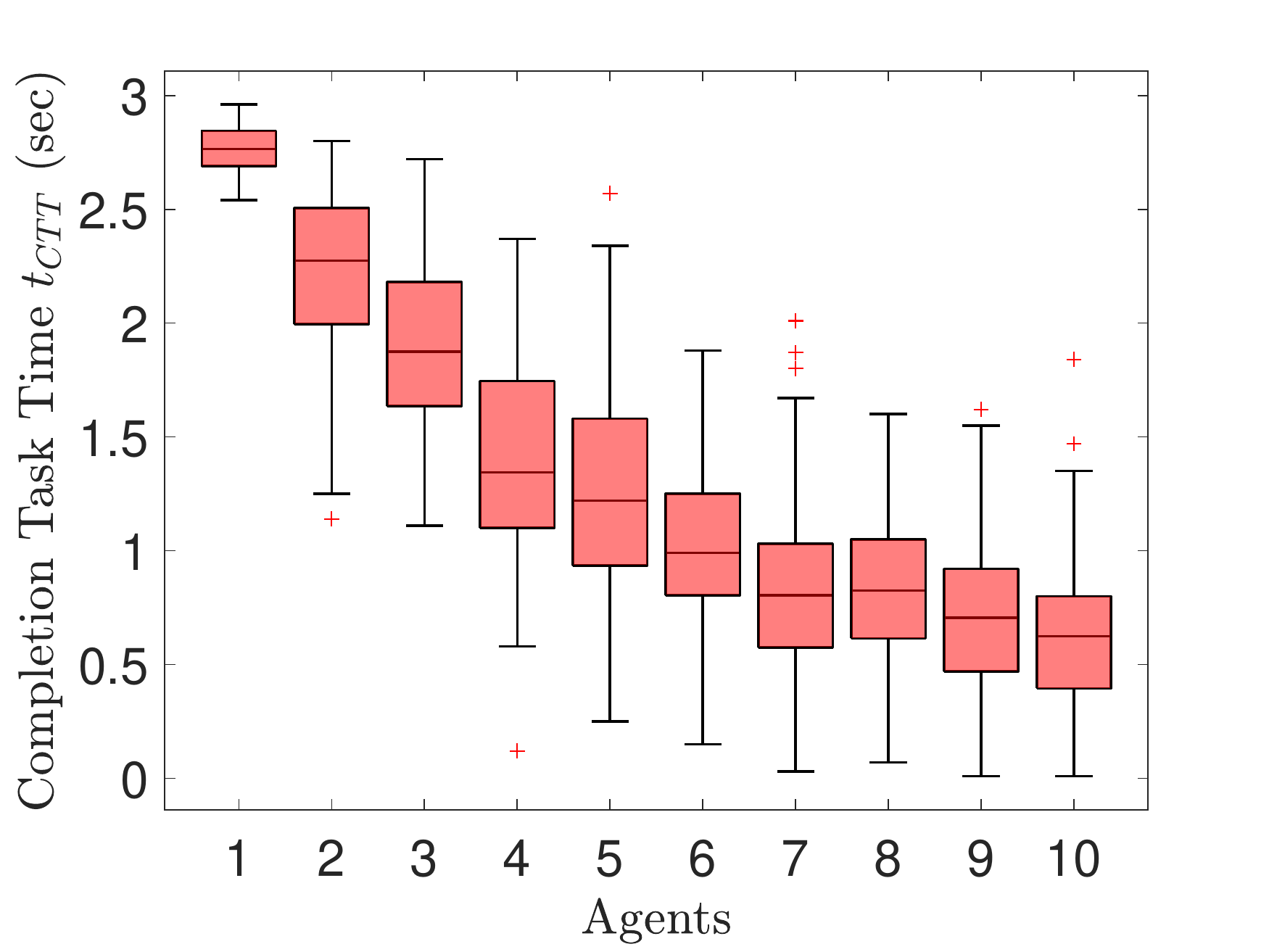}
         \caption{ }
         \label{fig:bp1b}
     \end{subfigure}
   %  }
     \caption{Boxplot statistics of completion task time $t_{CTT}$ versus number of agents for Volcano (a), Archipelago (b) distributions.}
     \label{bp1}
 \end{figure*}
 \begin{figure*}[t]
     \centering
   %  \framebox{
     \begin{subfigure}{0.45\textwidth}
         \centering
         \includegraphics[trim={0cm 0 0cm 0},clip, scale =0.38]{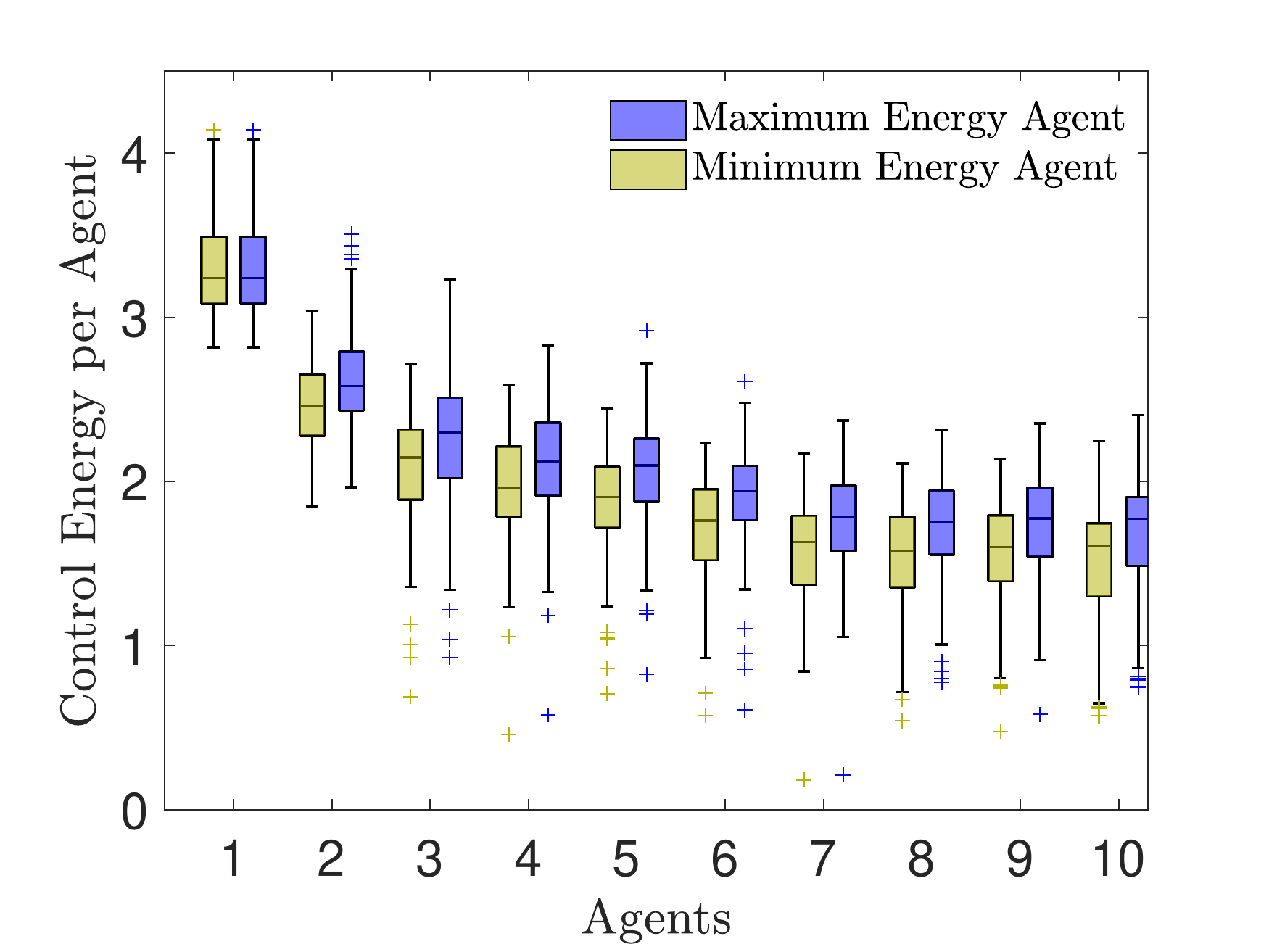}
         \caption{ }
         \label{fig:bp2a}
     \end{subfigure}
     \hfill
     \begin{subfigure}{0.45\textwidth}
         \centering
         \includegraphics[trim={0cm 0 0cm 0},clip, scale =0.38]{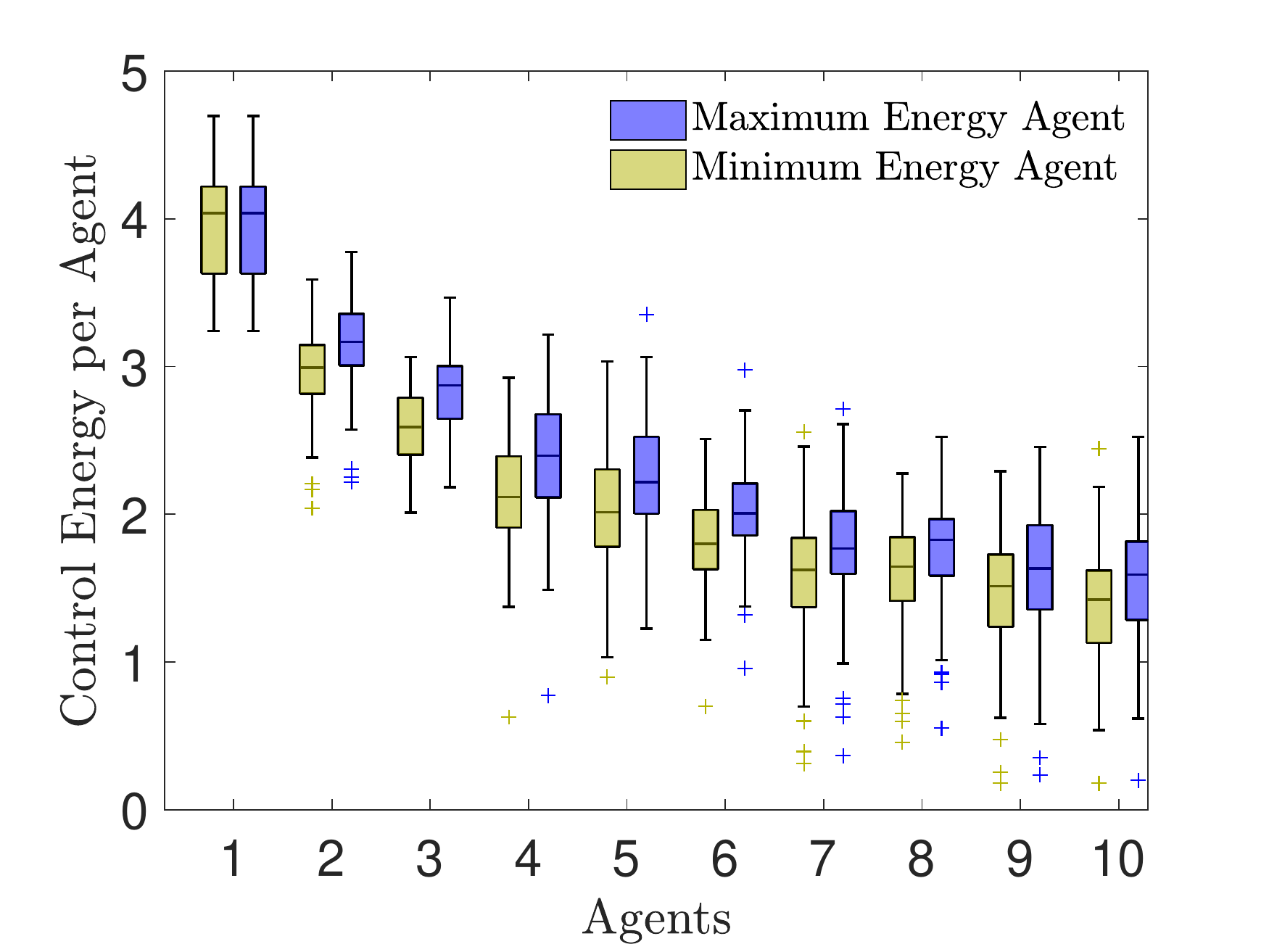}
         \caption{ }
         \label{fig:bp2b}
     \end{subfigure}
   %  }
     \caption{Boxplot statistics of control energy per agent versus number of agents for Volcano (a), Archipelago (b) distributions.}
     \label{bp2}
 \end{figure*}
  \begin{figure*}[t]
     \centering
     %\framebox{
     \begin{subfigure}{0.45\textwidth}
         \centering
         \includegraphics[trim={0cm 0 0cm 0},clip, scale =0.38]{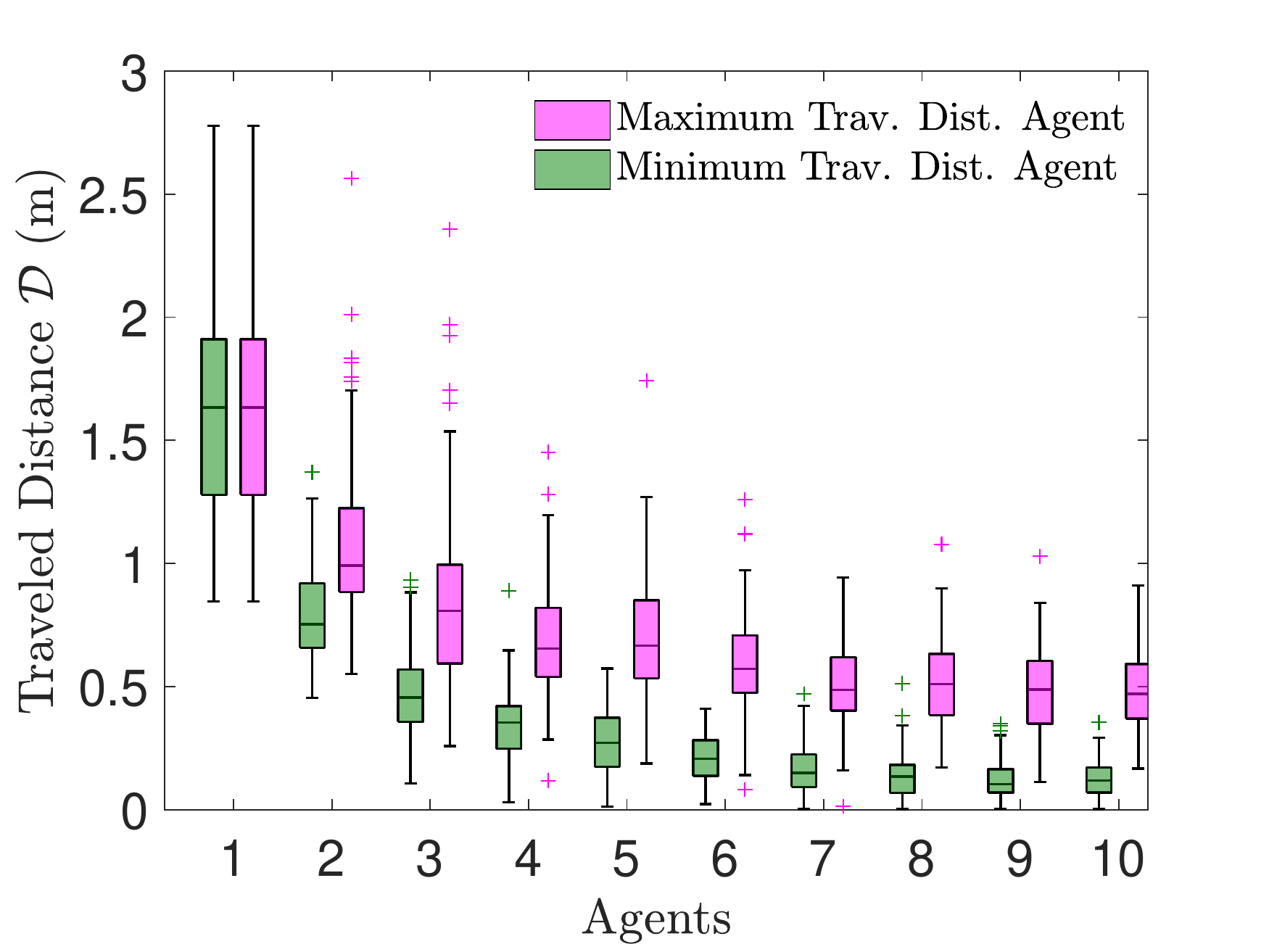}
         \caption{ }
         \label{fig:bp3a}
     \end{subfigure}
     \hfill
     \begin{subfigure}{0.45\textwidth}
         \centering
         \includegraphics[trim={0cm 0 0cm 0},clip, scale =0.38]{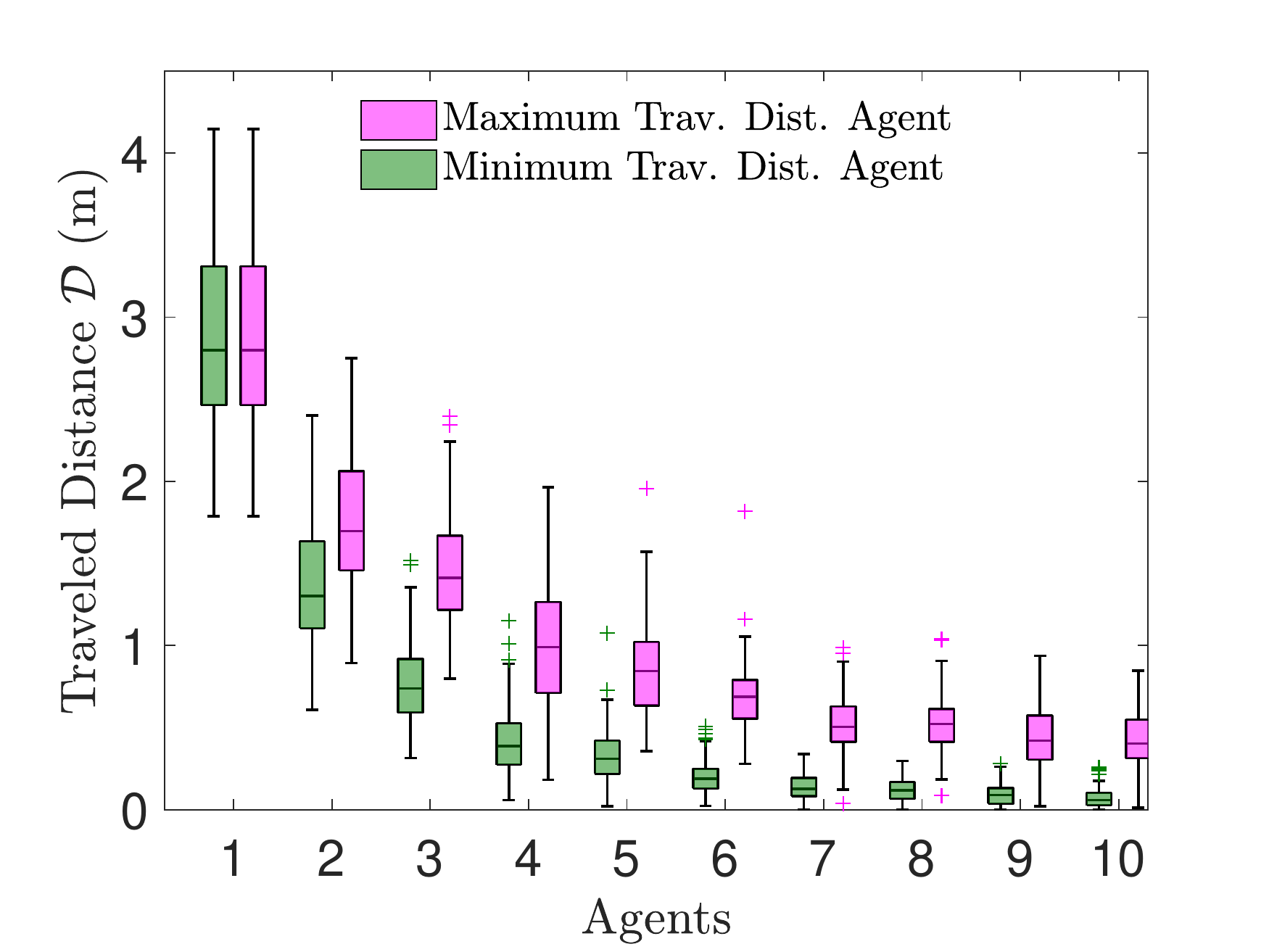}
         \caption{ }
         \label{fig:bp3b}
     \end{subfigure}
     %}
     \caption{Boxplot statistics of traveled distance $\mathcal{D}$ versus  number of agents for Volcano (a), Archipelago (b) distributions.}
     \label{fig:bp3}
 \end{figure*}

Figure \ref{fig:isl} presents similar analyses for the Archipelago
distribution. The single-agent visits all four modes spending more time in regions close to the peaks of the modes (Figure \ref{fig:isl_1}), achieving a completion time ($t_{CTT} = 3.17 \text{ sec}$) which is again very close to the horizon $T$ of the mission.
% Again it requires a high percentage of time in order to consider an efficient exploration as  $t_{CTT} = 3.25 \text{ sec}$. As far as the multi-agent system is
Collaborative planning is again observed for the 5-agent case in Fig. \ref{fig:isl_2}, where it is noted that the five initial trajectories are chosen very close to each other. Nonetheless, the optimized trajectories are distinct and allow the swarm to efficiently cover all modes ( $t_{CTT} = 1.65 \text{ sec}$).
%Meanwhile, two agents  visit their corresponding modes more than once. Consequently, the outcome is
% an efficient search domain that requires a relatively small fraction of time to accomplish the exploration task as shown by
% $t_{CTT} = 1.65 \text{ sec}$  in Fig. \ref{fig:isl_2}.
It is again instructive to observe the trend exhibited by the optimal ergodic metric $\mathcal{E}_{opt}$ (Fig. \ref{fig:isl_3}). In a qualitatively similar manner as before, the single-agent trajectory %shows first a good trend in the profile of $\mathcal{E}_{opt}$
initially a decrease in ergodicity
due to the coverage of the mode on the top-right. However, when it is directed to the second it spends a great amount of time in the intermediate low-valued information region, determining an over-shoot in the plot. This is clearly avoided in the multi-agent scenario configuration.

%Also, note that in Fig. \ref{fig:vol_3} and
% Fig. \ref{fig:isl_3}, we omit the $\mathcal{E}_{opt}(t)$ values at $t=0$ for clarity of presentation. That is  $\mathcal{E}_{opt}(t=0) = 5.307$ for the Volcano
% distribution and $\mathcal{E}_{opt}(t=0) = 4.878$ for the Archipelago distribution, while it is common by definition for any number of agents and depends on
% distribution characteristics.

\begin{figure}[ht]
    \centering
   \begin{subfigure}{.5\columnwidth}
       \centering
      \includegraphics[trim={0cm 0 0cm 0},clip, scale =0.3]{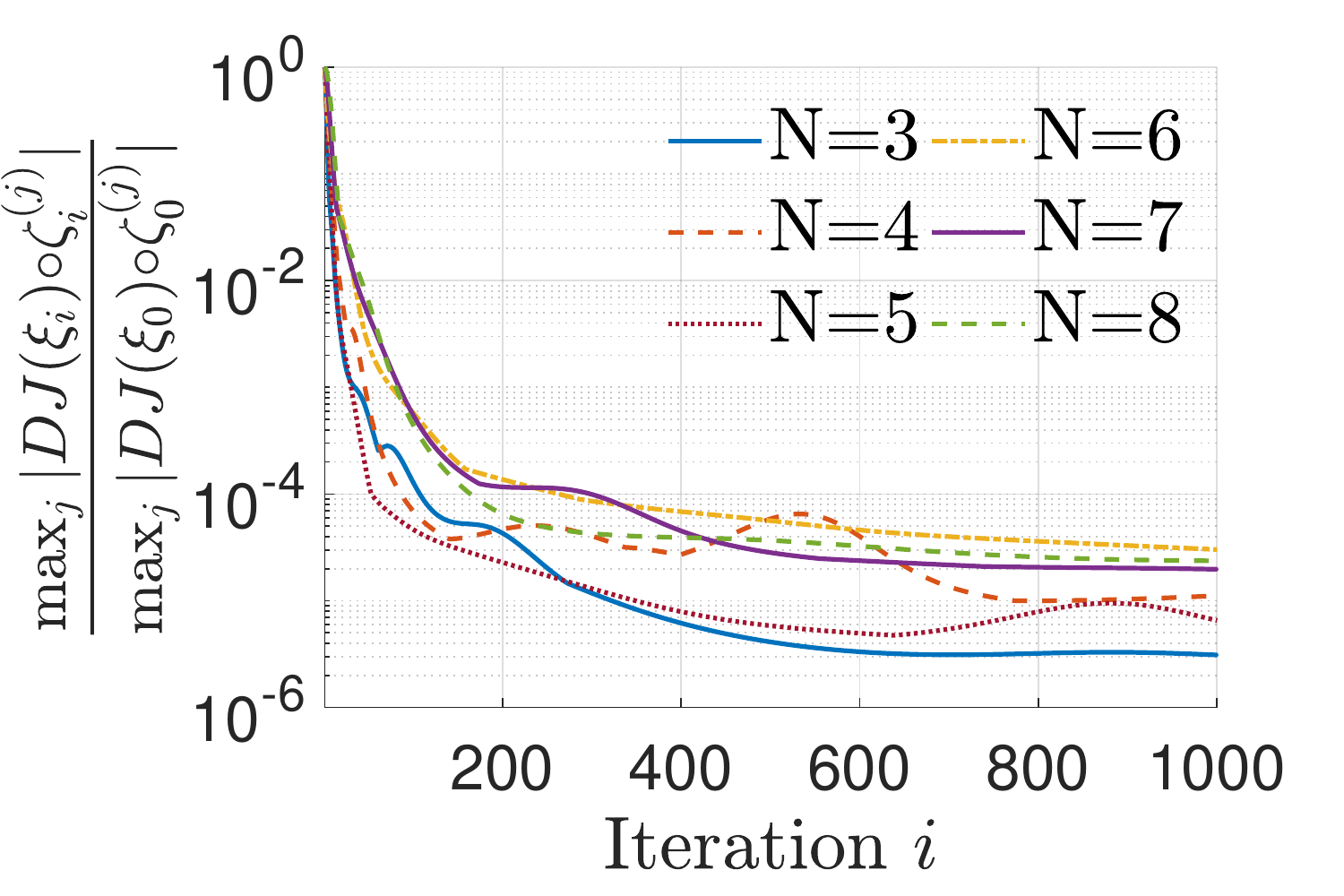}
      \caption{}
      \label{fig:cvra}
   \end{subfigure}%
    \hfill
  \begin{subfigure}{.5\columnwidth}
      \centering
      \includegraphics[trim={0cm 0 0cm 0},clip, scale = 0.3]{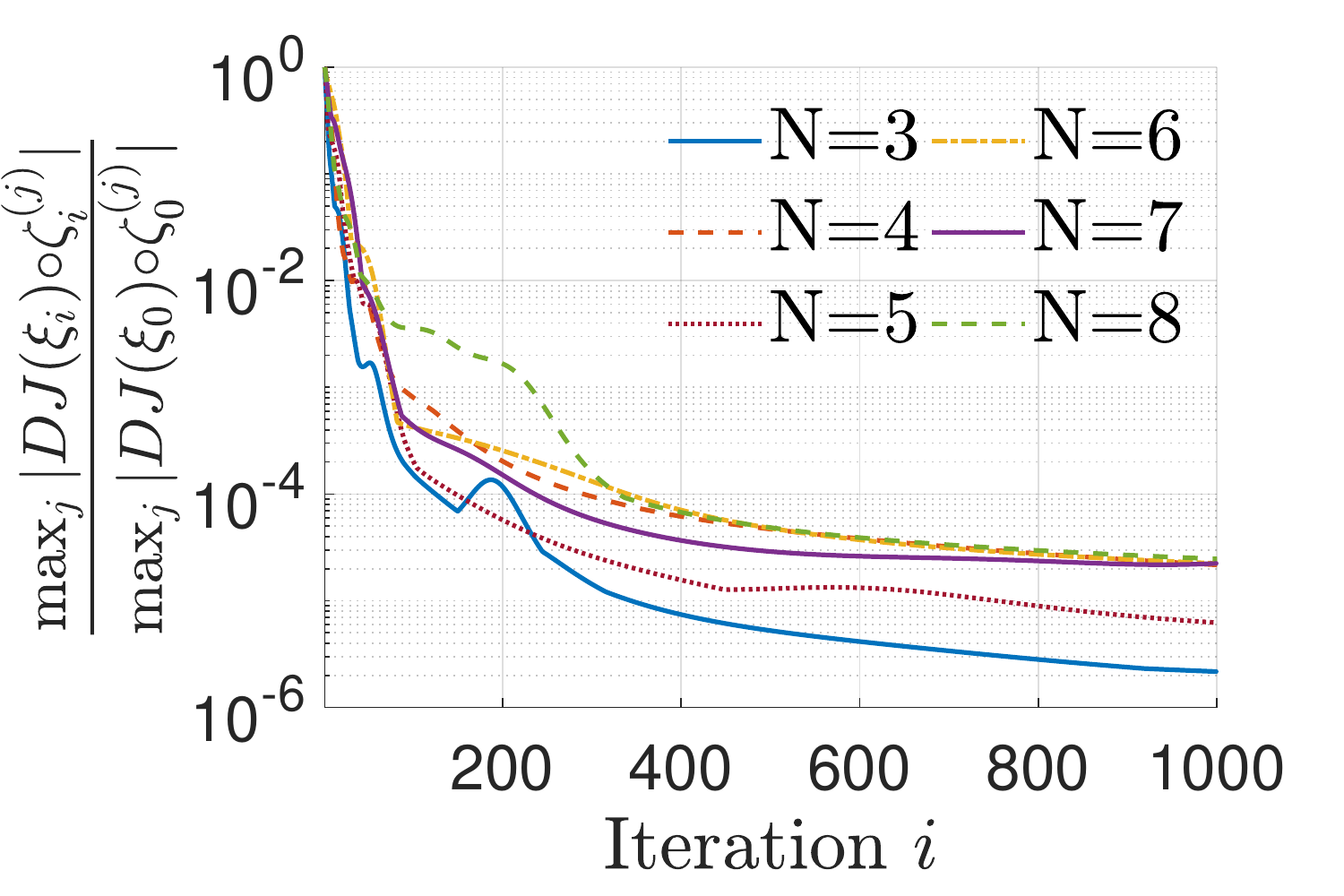}
      \caption{}
      \label{fig:cvrb}
  \end{subfigure}%
\caption{Convergence of the directional derivative versus iteration rounds for Volcano (a) and Archipelago (b) distributions.}
\label{fig:cvr}
\end{figure}

%Monte Carlo simulation
\subsection{Quantitative aspects and initialization effects}
%Although previous case studies depict insightful aspects of the methods, it is important to provide a systematic approach for the algorithmic evaluation.
%Therefore, the proposed algorithms are extensively assessed through Monte Carlo randomized initializations.
%One important aspect of the proposed trajectory planning approach is that, due to the non-convexity of the optimization problem, the final result will be sensitive to the initial trajectory used in the iterative.
%For this reason,
We provide here a comprehensive analysis of the effect of increasing the number of agents and of choosing the initial conditions of the agents' state on three performance metrics.
The number of considered agents is varied from 1 to 10, and, in each case, the optimization algorithm is run using 100 random initializations for the agents initial location $x_0$. Precisely, we sample each initial condition from a uniform distribution such that $\theta (t=0) \in [0, 2\pi]$, $X(t=0) \in [0.05, 0.95]$ and $Y(t=0) \in [0.05, 0.95]$. It is noted that this also randomizes the feasible trajectory used to initialize the optimization. As mentioned earlier, the trajectories are circles with given radius, and the selection of a point $x_0$, together with the constraint that the circular trajectory is feasible for the dynamics, uniquely determines its location in the field. These analyses thus shed also some light on the effect of the initial trajectory on the final result, which is an important aspect due to the non-convexity of the optimization problem.

Results are shown in terms of three performance metrics, namely the completion task time
$t_{CTT}$, the control energy per agent, and the traveled distance $\mathcal{D}$. %Overall, we employ 100 different random initializations per case and, thus, 1000 in total.
%The employed parameters are obtained through
% Table \ref{tab:tuning_parameters} values, while we use $r = 3$ for distance cost weight.  In order to speed up the process, we use a
% convergence tolerance value $\epsilon_2 = 95 \%$, while the completion task time $t_{CTT}$ is defined based on $\epsilon_{opt} = 98.5 \%$.
The \Resp{inter-agent distance} penalty term  and optimal temporal ergodicity reduction tolerance are set to $r = 3$ and $\epsilon_{opt} = 98.5$, respectively.

Figures \ref{fig:bp1a} and \ref{fig:bp1b} show box plots of the completion task time statistics against number of agents for Volcano and Archipelago distributions,
respectively. As expected, increasing the number of agents reduces the completion task time leading to more efficient exploration schemes.
In the Volcano case, above 8 agents no further decrease is observed, suggesting that there is a distribution-dependent threshold for the largest number of agents giving an advantage in completion time. \Resp{It is noted that $99.7\%$  and $92.4\%$ of the cases in Volcano and Archipelago distributions, respectively, have achieved an ergodicity reduction $\mathcal{E}_r$ above $95\%$.}

%In terms of control energy, we illustrate  the box plot statistics of control energy per agent in
Figure \ref{bp2} presents the control energy performance for the most and least energy consuming agents.% Volcano and Archipelago distributions, respectively. Actually, the statistics of the most energy consuming agent and the least energy consuming agent are
%presented.
 This analysis highlights the advantageous distribution of energy for a multi-agent system. Indeed, by increasing the agents number,
each individual agent consumes lower energy and, thus, the energy pool is distributed  efficiently among the agents.

Figure \ref{fig:bp3} finally shows the traveled distance box plot statistics. Recall that this metric is not explicitly targeted in the optimization, but it is nonetheless of practical interest to monitor it.
%an implicit performance metric.
As before, since this metric is a function of the agent, the minimum and maximum traveled distances are presented.Increasing the number of agents clearly leads to a decrease in the traveled distance. It is worth noting that it also leads to a decrease in the dispersion of this metric, both in terms of number of whiskers and box plots width. %Note finally that the traveled distance values of Volcano distributions are lower compared to Archipelago values as  Volcano density  is centralized into one peak mode.

\Resp{
\subsection{Convergence study}
Figure \ref{fig:cvr} reports the results of an investigation of the convergence properties of the decentralized algorithm.
%The directional derivative $D_{ \xi_{i}^{(j)}}  J\left(\xi_{i}\right) \circ \zeta_i ^{(j)}$ is taken as measure of optimality of the solution at iteration $i$.
Specifically, the largest absolute value of the directional derivative across agents, i.e.
$\underset{j \in \{1, \cdots, N\}} {\max }|D_{ \xi_{i}^{(j)}}  J\left(\xi_{i}\right) \circ \zeta_i ^{(j)}|$, normalized by the initial trajectory directional derivative at $i=0$, is presented as a metric of optimality.
%The directional derivative is computed using the true agent trajectories $\xi^{(j)}$.
%and not the estimated ones, in order to observe the true convergence property of the algorithm.
The analyses are done on the graphs in Figure \ref{fig:erdos-renyi} for a randomly generated initial trajectory and initial condition $x_0$. %  $2\leq N\ \leq 8$ agents and randomly produced initial locations $x_0$.

A sublinear rate can be recognized in both examples. It is also observed that, probably due to the heuristic decentralization of the Armijo line search, there is no monotonic decrease in the optimality metric. However, this does not compromise convergence. Analyses have also been carried out for two extreme network topologies, namely complete and line graphs, which showcased qualitatively similar trends (data not shown).
}

\section{CONCLUSION}\label{sVI}
The paper presents a new approach to design trajectories of multi-agent systems for ergodic exploration of \Resp{stationary target distributions}. To this aim, an objective function comprising three distinct terms is defined, and a decentralized optimization algorithm is proposed to minimize it. Two examples of distributions are considered in numerical experiments, and results are shown to demonstrate the validity of the approach and support the advantages of the proposed solution.
%coupled with steepest descent optimization
%and projection operation. Multiple agents collaborate to explore efficiently search domains and outperform single-agent exploration.
%Monte Carlo simulations illustrate the advantages of multi-agent systems with respect to a variety of performance metrics.
%We demonstrate that by several metrics,
The multi-agent algorithm enables more efficient exploration strategies compared to the single-agent case. Importantly, the shared ergodic metric allows multiple agents to explore cooperatively in order to search efficiently the domain with a low completion task time and a more efficient use of energy.
%\Resp{The method is applicable to stationary target distributions, and future work shall investigate its extension to scenarios where measurements can be used online to update the exploration strategies. Experimental validation of the performance of the algorithm in a real-world environment is another important follow-up activity.}
\Resp{Future work shall investigate an extension of the proposed method to scenarios where the target distribution can be updated online with measured data, as well as an experimental validation in a real-world environment.}

\addtolength{\textheight}{-12cm}   % This command serves to balance the column lengths

%%%%%%%%%%%%%%%%%%%%%%%%%%%%%%%%%%%%%%%%%%%%%%%%%%%%%%%%%%%%%%%%%%%%%%%%%%%%%%%%

%%%%%%%%%%%%%%%%%%%%%%%%%%%%%%%%%%%%%%%%%%%%%%%%%%%%%%%%%%%%%%%%%%%%%%%%%%%%%%%%

%%%%%%%%%%%%%%%%%%%%%%%%%%%%%%%%%%%%%%%%%%%%%%%%%%%%%%%%%%%%%%%%%%%%%%%%%%%%%%%%
%\section*{APPENDIX}

%\section*{ACKNOWLEDGMENT}

%%%%%%%%%%%%%%%%%%%%%%%%%%%%%%%%%%%%%%%%%%%%%%%%%%%%%%%%%%%%%%%%%%%%%%%%%%%%%%%%


\end{document}